\documentclass[MPS,STIX1COL]{WileyNJD-v2}

\articletype{Article Type}%

\received{26 April 2016}
\revised{6 June 2016}
\accepted{6 June 2016}

\usepackage[caption=false]{subfig}
\usepackage{mathtools}
\usepackage[bbgreekl]{mathbbol}
\usepackage{tikz}
\usetikzlibrary{arrows.meta}

\def\R{\mathbb{R}}
\def\Rds{\mathbb{R}^{d \times d}_{\sym}}

\def\fp{ : }
\def\vec#1{\boldsymbol{#1}}
\def\bu{\vec{u}}
\def\bv{\vec{v}}
\def\bw{\vec{w}}
\def\bn{\vec{n}}

\def\tens#1{\pmb{\mathsf{#1}}}
\def\tr{\mathop{\mathrm{tr}}\nolimits}
\def\sym{\mathop{\mathrm{sym}}\nolimits}
\def\BD{\tens{D}}
\def\Bsigma{\tens{\sigma}}

\def\linstress{\bbtau}
\def\linstrain{\bbespilon}
\def\identity{\mathbb{I}}
\def\cA{\mathcal{A}}

\def\diver{\mathop{\mathrm{div}}\nolimits} 
\def\dotd#1{\dot{\overline{#1}}}
\def\dtau{\mathrm{d}_t^{\tau}}

\newcommand{\ymod}{\mathrm{E}\,}

\newcommand{\dd}[2]{\frac{{\rm d} #1}{{\rm d} #2}} 
\newcommand{\Heaviside}{H} 
\newcommand{\Young}{\ymod} 
\newcommand{\yield}{y} 
\newcommand{\absnorm}[1]{\left| #1\right|} 

\begin{document}

\title{Numerical approximation of a thermodynamically complete rate-type model for the elastic--perfectly plastic response}

\author[1,2]{Pablo Alexei Gazca-Orozco*}

\author[2]{V\'it Pr\r{u}\v{s}a}

\author[2]{Karel T\r{u}ma}

\authormark{Pablo Alexei Gazca-Orozco \textsc{et al}}

\address[1]{\orgdiv{Department of Mathematics}, \orgname{University of Freiburg}, \orgaddress{\state{Ernst-Zermelo-Stra\ss e, 79104, Freiburg}, \country{Germany}}}

\address[2]{\orgdiv{Charles University}, \orgname{Faculty of Mathematics and Physics}, \orgaddress{\state{Sokolovsk\'a 83, Praha, CZ 186~75}, \country{Czech Republic}}}

\corres{*Pablo Alexei Gazca--Orozco, Department of Mathematics, University of Freiburg, Ernst-Zermelo-Stra\ss e, 79104, Freiburg, Germany. \email{alexei.gazca@mathematik.uni-freiburg.de}}


\abstract[Abstract]{We analyse a numerical scheme for a system arising from a novel description of the standard elastic--perfectly plastic response. The elastic--perfectly plastic response is described via rate-type equations that do not make use of the standard elastic-plastic decomposition, and the model does not require the use of variational inequalities. Furthermore, the model naturally includes the evolution equation for temperature. We present a low order discretisation based on the finite element method. Under certain restrictions on the mesh we subsequently prove the existence of discrete solutions, and we discuss the stability properties of the numerical scheme. The analysis is supplemented with computational examples.
}

\keywords{Rate-type constitutive relations, perfect plasticity, finite element method, thermodynamically consistent models}

\maketitle


\section{Introduction}

A rate-independent hysteretic response is frequently encountered in various engineering applications such as electrical engineering, geomechanics and mechanical engineering. Each of these research communities have developed its own approaches to the modelling of the hysteretic response, see~\cite{visintin.a:differential}, \cite{ismail.m.ikhouane.f.ea:hysteresis} and \cite{pei.j.gay-balmaz.f.ea:connecting} for a list of various hysteretic models and a discussion of their relations. In solid mechanics the prime example of a rate-independent hysteretic response is the elastic-plastic response, see for example~\cite{bruhns.ot:history} and~\cite{jarecki.d.srinivasa.ar.ea:rate-independent} for comments on the historical development of plasticity theory. In the present contribution we work with a novel model for the standard elastic--plastic response, see~\cite{rajagopal.kr.srinivasa.ar:inelastic,rajagopal.kr.srinivasa.ar:implicit*1}  and~\cite{CP.2020}, and we focus on mathematical aspects of the model. In particular we prove solvability of the corresponding spatially discretised system of governing partial differential equations.

Before we proceed with the numerical analysis, let us briefly comment on the status of the considered model. Concerning the elastic--plastic response of metals, the predominant modelling approach is based on the concepts of the \emph{elastic--plastic decomposition}, the \emph{flow rule} and the \emph{yield condition}, which in turn leads to a characterisation of the elastic--plastic response using the concepts of optimisation theory, see especially~\cite{SH.2006} and for further discussion also~\cite{scalet.g.auricchio.f:computational}. Concerning the elastic-plastic response of non-metallic materials such as soils, the situation is different, see, for example, \cite[Ch.~8]{hashiguchi.k:elastoplasticity} for a relatively recent critical review of some popular models and \cite{hashiguchi.k:mechanical} for an early discussion of the same. These materials typically exhibit ``diffuse yielding behaviour''/ ``smooth elastic-plastic transition'', see~\cite{jarecki.d.srinivasa.ar.ea:rate-independent} and the discussion therein,  which means that the transition from the elastic to the plastic regime is not sharp, but it progresses gradually, hence the concept of \emph{sharp yield condition} must be abandoned. (Some exotic alloys however seem to exhibit the same behaviour as well, see~\cite{besse.m.castany.p.ea:mechanisms}, \cite{mozafari.f.thamburaja.p.ea:rate} or~\cite{saito.t.furuta.t.ea:multifunctional}.)  In this case the elastic--plastic response is typically modelled using \emph{rate-type} equations designed in such a way that the whole model still predicts the \emph{rate-independent} behaviour.

We shall investigate the family of models introduced in~\cite{rajagopal.kr.srinivasa.ar:inelastic}. This class of models belongs to the class of rate-type models, but it goes one step further. It also abandons the concept of elastic--plastic decomposition. In particular, the models do not use the traditional concept of 
strain decomposition to the elastic and plastic part, see~\cite{sadik.s.yavari.a:on} or~\cite{steigmann.dj:primer} and references therein; the models stemming from~\cite{rajagopal.kr.srinivasa.ar:inelastic} work with the stress and the strain only. In a one-dimensional setting the stress--strain relation is given by the rate-type equation
\begin{equation}
  \label{eq:1}
  \dd{\sigma}{t}
  =
  \Young
  \left[
    1
    -
    \Heaviside
    \left(\sigma \dd{\varepsilon}{t} \right)
    \Heaviside
    \left(
      \absnorm{\sigma} - \sigma_{\yield}
    \right)
  \right]
  \dd{\varepsilon}{t}
  ,
\end{equation}
where $\sigma$ denotes the stress, $\varepsilon$ denotes the strain, $\sigma_{\yield}$ denotes the yield stress, $\Young$ denotes the Young modulus and $\Heaviside$ denotes the Heaviside step function~\eqref{eq:Heaviside}. (Compare with the standard  models that lead to optimisation problems, see~\cite[Ch.~1]{SH.2006}.) The rate-type stress--strain relation~\eqref{eq:1} is clearly rate-independent, and it leads to the standard elastic--perfectly plastic response. Indeed, if the yield stress is reached, $\absnorm{\sigma} = \sigma_{\yield}$, and if the material is being loaded, $\sigma \dd{\varepsilon}{t} \geq 0$, then~\eqref{eq:1} reduces to
\begin{equation}
  \label{eq:2}
  \dd{\sigma}{t} = 0,
\end{equation}
hence the stress $\sigma$ remains constant and equal to the yield stress value. This is the plastic flow regime. On the other hand, if the stress is below the yield stress value, $\absnorm{\sigma} < \sigma_{\yield}$, \emph{or} if the material is being unloaded, $\sigma \dd{\varepsilon}{t} < 0$, then~\eqref{eq:1} reduces to
\begin{equation}
  \label{eq:3}
  \dd{\sigma}{t} = \Young \dd{\varepsilon}{t}.
\end{equation}
This is the standard elastic response rewritten in terms of rates. Indeed, equation~\eqref{eq:3} is just the time derivative of Hooke law $\sigma = \Young \varepsilon$. Once we have~\eqref{eq:2} and \eqref{eq:3}, it is straightforward to see that the cyclic change of strain leads to the standard hysteretic behaviour.

An important feature of the model~\eqref{eq:1} is that the second Heaviside function
$
\Heaviside
\left(
  \absnorm{\sigma} - \sigma_{\yield}
\right)
$ can be replaced by a smoothed version thereof, which allows one to easily deal with the ``diffuse yielding behaviour'', see~\cite{rajagopal.kr.srinivasa.ar:inelastic} and for further comments also~\cite{jarecki.d.srinivasa.ar.ea:rate-independent}. Furthermore, the family of one-dimensional models based on the rate-type equation~\eqref{eq:1} can be extended to the fully three-dimensional finite deformations setting, see~\cite{rajagopal.kr.srinivasa.ar:implicit*1} and~\cite{CP.2020}. 

Finally, the finite deformation version of the models can be shown to be thermodynamically consistent, see~\cite{CP.2020}. This implies that the energy conversions in the material are fully characterised. In particular, the heat generated in the inelastic processes is known, and the models allow one to study fully coupled thermomechanical processes in the finite strain setting. Despite their importance, such coupled thermomechanical processes are rarely studied, especially in the case of rate-type models for soils, see~\cite{hashiguchi.k:nonlinear,hashiguchi.k:foundations,hashiguchi.k:elastoplasticity}, and the situation is only slightly better for metals, see~\cite{rosakis.p.rosakis.aj.ea:thermodynamic} for an early example thereof.

In the present work, we focus on a model of type~\eqref{eq:1} that arises in the \emph{small strain approximation} of a finite deformation model based on~\eqref{eq:1}, while we study both withe ``sharp yielding behaviour'' as well as ``diffuse yielding behaviour''.  The model is described in~\cite{CP.2020}, and it focuses on the core features of elastic--plastic material response. It is a relatively simple model without additional features such as kinematic/isotropic hardening, and as such it is suitable for a proof-of-concept numerical analysis of this class of models. In particular, we prove solvability of the equations arising from the spatial discretisation of the corresponding partial differential equations, and we also investigate stability properties of the corresponding numerical scheme.

\section{Model description}
Let the computational domain $\Omega$---which is tantamount to the reference stress-free configuration of the body of interest---be an open bounded subset of $\R^d$, with $d\in\{2,3\}$, whose boundary $\partial\Omega$ is Lipschitz. Concerning the boundary conditions for the mechanical quantities, we assume that the domain is disjointly divided into a Dirichlet (displacement) $\partial\Omega_D$ and a Neumann (traction) $\partial\Omega_N$ component. For the thermal part we prescribe the zero Neumann (no-flux) boundary condition everywhere on the boundary~$\partial\Omega$; see below for details.  

As shown in \cite{CP.2020} the standard elastic--perfectly plastic response with \emph{von Mises} yield criterion in the small strain regime can be described by the following system of equations posed on the space-time domain $Q:= (0,T)\times \Omega$:

\begin{subequations}\label{eq:pde_Tv}
\begin{align}
    \rho_\star \dotd{\bv} &-  \diver\linstress 
    = \rho_\star\mathbf{f} \qquad \quad &&\text{ in }(0,T)\times \Omega,  \label{eq:pde_momentum}\\
    \frac{1}{\ymod}( (1+\nu) \dotd{\linstress} - \nu (\tr\dotd{\linstress})\mathbb{I}) = \dotd{\linstrain}
                          &- H(\linstress\fp \dotd{\linstrain} )H(|\linstress_\delta|^2 - \kappa_\star^2)\dotd{\linstrain}
                          &&\text{ in }(0,T)\times\Omega, \label{eq:pde_CR}
    \end{align}
\end{subequations}
plus initial conditions $\bv(0,\cdot)=\bv_0(\cdot)$ and $\linstress(0,\cdot)=\linstress_0$, and boundary conditions $\bv|_{\partial\Omega_D}=\bv_b$ and $\linstress\bn|_{\partial\Omega_N}=\mathbf{t}_b$, where $\overline{(\partial\Omega_D)}\cup \overline{\partial\Omega_N}=\partial\Omega$.

Here $\bu$ denotes the displacement, $\linstrain = \linstrain(\bu) := \frac{1}{2}(\nabla \bu + \nabla \bu^\top)$ denotes the linearised strain operator (infinitesimal strain tensor), $\linstress$ denotes the stress tensor, and $\mathbb{A}_\delta := \mathbb{A} - \frac{1}{d}\tr(\mathbb{A})\mathbb{I}$ denotes the traceless part of the corresponding tensor $\mathbb{A}$. The symbol $\left( \mathbb{A}\fp\mathbb{B} \right):= \tr \left({\mathbb{A}\mathbb{B}^{\mathrm{T}}} \right)$ denotes the matrix scalar product. The function $H$ is the classical Heaviside function, defined as
\begin{equation}\label{eq:Heaviside}
H(s):= \left\{
\begin{array}{lc}
  1, & s\geq 0,\\
  0, & s<0.
\end{array}\right.
\end{equation}
All the quantities of interest are  functions of the position in the reference configuration $\mathbf{X}\in \Omega$ and time $t\in (0,T)$; the dot represents the time derivative
\begin{equation*}
\dotd{\mathbb{A}} := \frac{\partial}{\partial t}\mathbb{A}(t, \mathbf{X}).
\end{equation*}
The symbols $\ymod$, $\nu$, $\kappa_\star$ denote material parameters, namely Young modulus, Poisson ratio, and yield stress; the density in the reference configuration is denoted by $\rho_\star$ and we assume that $\rho_\star \geq \rho_\star^-$, for some positive constant $\rho_\star^-$.

The first equation \eqref{eq:pde_momentum} represents balance of momentum, and the second equation \eqref{eq:pde_CR} is the rate-type constitutive relation for the elastic--perfectly plastic response. Since in the small strain regime we have $\dotd{\linstrain} = \linstrain(\bv)$, we see that the system \eqref{eq:pde_Tv} is a system of evolution equations for the velocity field $\bv$ and the stress tensor $\linstress$.

The displacement $\bu$ and temperature $\theta$ can be computed post-hoc. Once the stress $\linstress$ and the velocity $\bv$ fields are known, it remains to solve
\begin{subequations}\label{eq:pde_ut}
\begin{align}
    \dotd{\bu} &= \bv  &&\text{ in }(0,T)\times\Omega,\\
    \rho_\star c_v \dotd{\theta} - \diver(\kappa_{\textrm{th}}\nabla \theta) = H(\linstress &\fp \dotd{\linstrain} )H(|\linstress_\delta|^2 - \kappa_\star^2) \linstress\fp \dotd{\linstrain}\  &&\text{ in }(0,T)\times\Omega \label{eq:pde_temperature},
  \end{align}
\end{subequations}
for $\bu$ and $\theta$. Here $c_v$ denotes the specific heat capacity at constant volume, and $\kappa_{\textrm{th}}$ is the thermal conductivity. The boundary conditions for the displacement are chosen to be consistent with those of $\bv$; i.e.\ if $\bu|_{\partial\Omega_D} = \bu_b$ then $\bv|_{\partial\Omega_D} = \bv_b := \dotd{\bu_b}$. For the temperature we impose the no-flux boundary condition, that is $\nabla\theta \cdot \bn|_{\partial\Omega}=0$. We note that the presented model is rather simple, the coupling between the governing equations for the thermal and mechanical quantities can be more involved---the yield stress can be a function of temperature, the elastic response can be designed in such a way that it captures the Gough--Joule effect, see \cite{gough.j:description}, \cite{joule.jp:on*2} and \cite{anand.l:constitutive} for a modern discussion thereof, and so forth. These more involved temperature related effects can be easily added into the model. (For example, the Gough--Joule effect can be captured by using the same Helmholtz free energy \emph{ansatz} as in~\cite{anand.l:constitutive}. See also~\cite{prusa.v.tuma.k:temperature} for a discussion in the case of a related rate-type model for a viscoleastic solid.) However, from the perspective of \emph{rigorous numerical analysis} these effects bring further complications, see also Remark~\ref{rmk:temp_dependent_parameters} at the end of Section~\ref{sec:discrete-formulation}. Since the more involved models for thermomechanical coupling might lead to difficulties in \emph{rigorous mathematical and numerical analysis}, and they would also require a careful discussion of possible (mathematically treatable) thermal effects, we deliberately ignore the more involved thermal effects. We focus only on the core analytical issues related to the novel rate-type model for the elastic--perfectly plastic response, and the simple temperature evolution equation~\eqref{eq:pde_temperature} is sufficient from this point of view.   

One of the challenging aspects of system \eqref{eq:pde_Tv} is the presence of the Heaviside function, since then one has to deal with a differential equation with a discontinuous nonlinearity. To alleviate this difficulty, we will employ a non-sharp yield condition. (In the terminology used in \cite{jarecki.d.srinivasa.ar.ea:rate-independent} this is tantamount to the ``diffusive yielding behaviour''.) The non-sharp yield condition means that the last term in \eqref{eq:pde_CR} is substituted by $H(\linstress\fp \dotd{\linstrain}) H_\epsilon(|\linstress_\delta|^2-\kappa_\star^2)\dotd{\linstrain}$, where $H_\epsilon$ is a regularised version of the Heaviside function. We consider three different options in this work, namely
\begin{subequations}
\begin{align}
  H^{(1)}_\epsilon(s) &:= \frac{1}{2} + \frac{1}{2}\frac{\frac{s}{\epsilon}}{\sqrt{1 + (\frac{s}{\epsilon})^2}}
                      &\qquad \epsilon>0,\ s\in \mathbb{R},\\
  H^{(2)}_\epsilon(s) &:= \frac{1}{2} + \frac{1}{2}\tanh\left(\frac{s}{\epsilon}\right)
                      &\qquad \epsilon>0,\ s\in \mathbb{R},\\
  H^{(3)}_\epsilon(s) &:= \frac{1}{2} + \frac{1}{\pi}\arctan\left(\frac{s}{\epsilon}\right)
                      &\qquad \epsilon>0,\ s\in \mathbb{R},
\end{align}
\end{subequations}
where $\epsilon$ is the regularisation parameter.

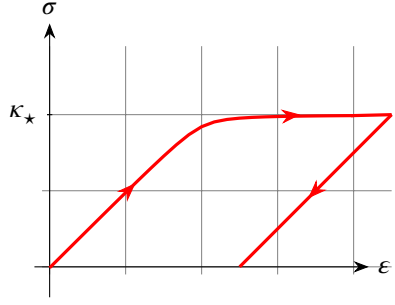
\begin{figure}
\centering
  \begin{tikzpicture}[domain=0:4] 
    \draw[very thin,color=gray] (-0.1,-0.1) grid (4.5,2.9);
    \draw[-{Stealth[scale=1.2]}] (-0.2,0) -- (4.2,0) node[right] {$\varepsilon$}; 
    \draw[-{Stealth[scale=1.2]}] (0,-0.4) -- (0,3.2) node[above] {$\sigma$};
    \draw[-{Stealth[scale=1.5]},color=red] (1.03,1.0) -- (1.13,1.1) ;
    \draw[-{Stealth[scale=1.5]},color=red] (3.5,1.0) -- (3.4,0.9) ;
    \draw[-{Stealth[scale=1.5]},color=red] (3.2,1.99) -- (3.3,2.0) ;
    \draw[color=red, very thick]    plot (\x,{2 + 0.5*(\x - 2 - sqrt((\x-2)*(\x-2) + 0.1))}) ;
    \draw[color=red, very thick] (3.965,1.987) -- (4.5,2);
    \draw[color=red, very thick] (4.5,2.0) -- (2.5,0);
    \draw (1pt, 2) -- (-1pt, 2) node[anchor=east,fill=white] {$\kappa_\star$};
  \end{tikzpicture}
  \caption{Non-sharp yield condition obtained as a consequence of a regularised Heaviside function $H_\varepsilon$.}
	\label{fig:non-sharp-yield}
\end{figure}

The qualitative one-dimensional behaviour during loading and unloading, for the stress $\sigma$ and strain $\epsilon$, that can be described by the non-sharp yield condition is depicted in Figure \ref{fig:non-sharp-yield}. (The magnitude of the regularisation parameter $\epsilon$ controls the ``sharpness'' of the corner on the loading curve.) We reiterate that the regularisation is in some physically relevant cases not artificial. In fact many materials exhibit such non-sharp yield conditions, see~\cite{jarecki.d.srinivasa.ar.ea:rate-independent} and the discussion therein. The freedom to model such non-sharp/diffusive yield condition is an advantage of the approach presented here, in contrast with the more widely used rate-independent models, where modelling non-sharp yield conditions is more cumbersome, see again~\cite{jarecki.d.srinivasa.ar.ea:rate-independent}.

\section{Discrete formulation}
\label{sec:discrete-formulation}
We employ the standard notation for Lebesgue spaces $(L^p(\Omega), \|\cdot\|_{L^p(\Omega)})$ and Sobolev spaces $(W^{1,p}(\Omega), \|\cdot\|_{W^{1,p}(\Omega)} )$. Let $\{\mathcal{T}_h\}_{h>0}$ be a family of shape-regular triangulations of $\Omega$ associated to a sequence of mesh sizes $h\to 0$; we assume here that $\Omega$ is a Lipschitz domain with polyhedral boundary, and also for simplicity we assume that the mesh is quasi-uniform, which implies that following inverse inequalites are available \cite{EG.2021.I},
\begin{subequations}\label{eq:inverse_inequalities}
\begin{gather}
\|\nabla\bv\|_{L^2(\Omega)} \leq 
c_{\textrm{inv}} h^{-1} \|\bv\|_{L^2(\Omega)}
\qquad \forall \bv\in V^h, \label{eq:inverse_ineq}\\
\|\bv\|_{L^2(\partial\Omega)} \leq c_{\textrm{tr}}h^{-1/2}\|\bv\|_{L^2(\Omega)}
\qquad \forall \bv\in V^h, \label{eq:inverse_ineq_bdry}
\end{gather}
\end{subequations}
where $c_{\textrm{inv}}, c_{\textrm{tr}}>0$ are positive constants independent of $h$; this quasi-uniformity assumption is not crucial, if desired one can apply instead local inverse inequalities.

The finite element spaces for the stress and velocity are chosen as 
\begin{gather*}
\Sigma^h = \{\Bsigma \in L^\infty(\Omega)^{d\times d}_{\sym} \, :\, \Bsigma|_K \in \mathbb{P}_0(K)^{d\times d}_{\sym},\, \forall K\in \mathcal{T}_h \} = \mathbb{DG}(0)^{d\times d}_{\sym},\\
V^h = \{\bv\in W^{1,\infty}(\Omega)^d \, :\, \bv|_K \in \mathbb{P}_1(K)^d ,\, \forall K\in \mathcal{T}_h ,\, \bv|_{\partial\Omega_D}=\mathbf{0}\},
\end{gather*}
that is, piecewise-linear Lagrange elements for the velocity and piecewise-constant approximations for the stress; here $\mathbb{P}_q(K)$ denotes the set of polynomials of degree at most $q$ on an element $K\in \mathcal{T}_h$. Since we are interested in approximating discontinuous terms, it is natural to employ lower order approximations, because higher degree polynomials could, in the absence of for example adaptivity, lead to unwanted oscillations.

For later use it is convenient to define the compliance operator $\mathcal{A}\colon \Rds\to \Rds$ as
\begin{equation}
\mathcal{A}(\bbsigma) := 
\frac{1}{\ymod}( (1+\nu) \bbsigma - \nu (\tr\bbsigma)\identity), \qquad \bbsigma\in \Rds.
\end{equation}
that is $\cA = \mathbb{C}^{-1}$, where $\mathbb{C}$ is the standard linear elasticity tensor. Since $\cA$ is positive definite, we can use it to define a norm on the space of discrete stresses $\Sigma^h$:
\begin{equation}
\|\bbsigma\|^2_{\cA} := \int_\Omega \cA\bbsigma \fp \bbsigma,
\qquad 
\bbsigma \in \Sigma^h.
\end{equation}
This norm is clearly equivalent to the $L^2$-norm; namely, $A_{\textrm{min}}\|\bbsigma\|^2_{L^2(\Omega)} \leq \|\bbsigma\|^2_{\cA} \leq A_{\textrm{max}}\|\bbsigma\|^2_{L^2(\Omega)}$, where $A_{\textrm{min}}$ and $A_{\textrm{max}}$ denote the minimum and maximum eigenvalues of $\cA$, respectively. Similarly, we equip the velocity space $V^h$ with the weighted norm
\begin{equation}
\|\bv\|^2_{\rho_\star} := \int_\Omega \rho_\star \bv\cdot \bv, \qquad \bv\in V^h.
\end{equation}
The corresponding weighted spaces of square integrable functions at the continuous level will be denoted by $L^2_{\rho_\star}(\Omega)$ and $L^2_{\cA}(\Omega)$.

Concerning the discretisation in time, we choose a time-step $\tau>0$, and we define a uniform time grid $t^k:= k\tau$, for $k\in \{1,\ldots, T/\tau\}$. (We can without loss of generality assume that $T/\tau\in \mathbb{N}$.) The system of governing equations is then discretised in time with the implicit Euler method; given a family of functions $\{\bv^k\}_{k\in\{0,\ldots, T/\tau\}}$ we define the discrete time derivative operator (or temporal difference quotient) as
\begin{equation}
\dtau \bv^k := \frac{\bv^k - \bv^{k-1}}{\tau},
\qquad
k\in \{1,\ldots,T/\tau\}.
\end{equation}

Now we are in the position to formulate a time-stepping scheme. We assume that the boundary datum $\bv_b$ can be seen as the restriction of some $\mathbb{CG}(1)$ function on $\Omega$, which we still denote by $\bv_b$, and we set $\linstress^0:= \linstress_0$ and $\bv^0:=\bv_0$. In the finite element formulation, assuming that approximations $\linstress^{k-1}\in \Sigma^h$ and $\bv^{k-1}\in  \bv_b + V^h$ at time $t^{k-1}$ have already been found, we look for $(\linstress^k_{h,\tau,\epsilon},\bv^k_{h,\tau,\epsilon}):=(\linstress^k,\bv^k)\in \Sigma^h\times (\bv_b+ V^h)$ such that
\begin{equation}\label{eq:discrete_Tv}
\begin{aligned}
\int_\Omega \cA(\dtau \linstress^k)\fp \bbsigma
- \int_\Omega \linstrain(\bv^k)\fp\bbsigma
+\int_\Omega &  H(\linstress^k \fp \linstrain(\bv^k) )H_\epsilon (|\linstress_\delta^k|^2 - \kappa_\star^2) \linstrain(\bv^k) \fp \bbsigma = 0   &&\forall \bbsigma \in \Sigma^h,\\
\int_\Omega\rho_\star \dtau \bv^k \cdot \bw 
               +\int_\Omega \linstress^k &\fp \linstrain(\bw) = \int_\Omega \rho_\star \mathbf{f}^k\cdot \bw
+ \int_{\partial\Omega_N} \mathbf{t}_b^k \cdot \bw
                             && \forall \bw\in V^h.
\end{aligned}
\end{equation}
Here $\mathbf{f}^k$ and $\mathbf{t}_b^k$ are approximations of $\mathbf{f}$ and $\mathbf{t}_b$ at time $t=t^k$, respectively. E.g.\ if $\mathbf{f}$ and $\mathbf{t}$ continuous, we can set $\mathbf{f}^k(\cdot):= \mathbf{f}(t^k,\cdot)$ and $\mathbf{t}_b^k(\cdot):= \mathbf{t}_b(t^k,\cdot)$.

The displacement and temperature problems \eqref{eq:pde_ut} are be discretised with piecewise linear Lagrange elements, that is the spaces of discrete displacements $U^h$ and discrete temperatures $\Theta^h$ are defined as
\begin{gather*}
U^h = \{\bu\in W^{1,\infty}(\Omega)^d \, :\, \bu|_K \in \mathbb{P}_1(K)^d ,\, \forall K\in \mathcal{T}_h,\, \bu|_{\partial\Omega_D}=\mathbf{0}\},\\
\Theta^h = \{\theta \in W^{1,\infty}(\Omega) \, :\, \theta|_K \in \mathbb{P}_1(K),\, \forall K\in \mathcal{T}_h \} = \mathbb{CG}(1).
\end{gather*}
In the discrete formulation we set $\bu^0:=\bu_0$ and $\theta^0:=\theta_0$, and for $k\in\{1,\ldots,T/\tau\}$, assuming that $\linstress^k\in \Sigma^h$, $\bv^k\in \bv_b + V^h$, $\bu^{k-1}\in \bu_b + U^h$ and $\theta^{k-1}\in\Theta^h$ are known, we look for $(\bu^k,\theta^k)\in (\bu_b + U^h)\times \Theta^h$ such that
\begin{equation}\label{eq:discrete_ut}
\begin{aligned}
\int_\Omega \dtau \bu^k \cdot \bw
- \int_\Omega \bv^k\cdot \bw
&= 0   &&\forall \bw \in U^h,\\
\int_\Omega \rho_\star c_v \dtau \theta^k \phi
               +\int_\Omega \kappa_{\textrm{th}}\nabla \theta^k\cdot \nabla\phi
= \int_\Omega H(\linstress^k \fp \linstrain(\bv^k) )H_\epsilon (&|\linstress_\delta^k|^2 - \kappa_\star^2) \linstress^k\fp \linstrain(\bv^k)\phi
                             && \forall \phi\in \Theta^h.
\end{aligned}
\end{equation}

\begin{remark}
Given the discontinuous nature of the stress space $\Sigma^h$, and noting that nonlinear functions of $\linstress^k$ remain piecewise constant, the equation for $\linstress^k$ in \eqref{eq:discrete_Tv} holds pointwise,
\begin{equation*}
\cA(\dtau\linstress^k) - \linstrain(\bv^k) + H(\linstress^k\fp \linstrain(\bv^k))H_\epsilon(|\linstress^k_\delta|^2 - \kappa^2_\star)\linstrain(\bv^k)  = 0 \qquad\text{in }\Omega.
\end{equation*}
This defines (implicitly) a mapping $\bv^k \to \tilde{\linstress}^k(\bv^k)$, which could be used to define a velocity-only problem
\begin{equation*}
\int_\Omega\rho_\star \dtau \bv^k \cdot \bw 
               +\int_\Omega \tilde{\linstress}^k(\bv^k) \fp \linstrain(\bw) = \int_\Omega \rho_\star \mathbf{f}^k\cdot \bw
+ \int_{\partial\Omega_N} \mathbf{t}^k \cdot \bw
                             \qquad \forall \bw\in V^h.
\end{equation*}
Using tools from automatic differentiation this can be solved, resulting in a strategy similar to the one employed traditionally, in which consistent tangents are employed in the linearisation \cite[Ch.\ 4.3.6]{SH.2006}. We do not pursue this further in this work.
\end{remark}

\begin{remark}
The discrete formulation \eqref{eq:discrete_Tv} involves $d(d+1)/2$ scalar $\mathbb{P}_0$ fields and $d$ scalar $\mathbb{P}_1$ fields; the problem  \eqref{eq:discrete_ut} involves in addition $d+1$ scalar $\mathbb{P}_1$ fields. Regarding the number of degrees of freedom, this is more computationally expensive than traditional formulations based on displacement and temperature. However, in addition to advantages from the modelling side, such as thermodynamic consistency, one advantage of the discrete formulation presented here is that the stress can be computed directly compared to traditional methods based on non-smooth constrained optimisation.
\end{remark}

\subsection{Existence of discrete solutions}
The goal in this section is to prove that numerical solutions to~\eqref{eq:discrete_Tv} exist. To help with this, we look first at the system in which both Heaviside functions are regularised. (The regularisation parameters are denoted as $\eta$ and $\epsilon$.) For simplicity we also assume that $\bv_b=\mathbf{0}$. Define the function $F_\eta \colon \Sigma^h\times V^h \to \Sigma^h \times V^h$ through the relation 
\begin{align}
\langle F_\eta(\linstress, \bv), (\bbsigma,\bw) \rangle 
& :=
 \int_\Omega  \cA(\linstress)\fp \bbsigma
- \tau \int_\Omega \linstrain(\bv)\fp\bbsigma
+  \tau \int_\Omega   H_\eta(\linstress \fp \linstrain(\bv) )H_\epsilon (|\linstress_\delta|^2 - \kappa_\star^2) \linstrain(\bv) \fp \bbsigma \notag\\
&\quad + \tau \int_\Omega \linstress\fp \linstrain(\bw)
+ \int_\Omega  \rho_\star\bv\cdot \bw
- \int_\Omega  \cA(\linstress^{k-1})\fp \bbsigma
- \int_\Omega \rho_\star \bv^{k-1}\cdot \bw \label{eq:discrete_residual} \\
&\quad - \tau \int_\Omega \rho_\star \mathbf{f}^k \cdot \bw 
-  \tau \int_{\partial\Omega_N} \mathbf{t}^k_b \cdot \bw. \notag
\end{align}
Note that the (regularised) discrete formulation can be written simply as $F_\eta(\linstress,\bv)=0$; note also that the function $F_\eta$ is continuous. If we manage to find a positive number $\hat{c}$, such that $\langle F_\eta(\linstress, \bv), (\linstress,\bv)\rangle \geq 0$, for all $(\linstress,\bv)\in \Sigma^h \times V^h$ with $\|\linstress\|^2_{\cA} + \|\bv\|^2_{\rho_\star} = \hat{c}$, then a corollary of Brouwer's fixed point theorem will guarantee the existence of a discrete solution \cite[Ch.\ 4, Cor.\ 1.1]{GR.1986}.
To this end, we take $(\bbsigma,\bw)=(\linstress,\bv)$ in the definition of $F_\eta$; this yields:
\begingroup
\allowdisplaybreaks
\begin{align*}
\langle F_\eta(\linstress, \bv), (\bbsigma,\bw) \rangle 
&\geq  
\int_\Omega \cA(\linstress)\fp \linstress 
+ \int_\Omega \rho_\star |\bv|^2
- \|\linstress^{k-1}\|_{\cA}\|\linstress\|_{\cA}
- \|\bv^{k-1}\|_{\rho_\star}\|\bv\|_{\rho_\star}
- \tau \|\mathbf{f}^{k}\|_{\rho_\star}\|\bv\|_{\rho_\star} \\
 &\quad- \frac{\tau c_{\textrm{tr}}}{(h\rho_\star^-)^{1/2}} \|\mathbf{t}^{k}_b\|_{L^2(\partial\Omega_N)}\|\bv\|_{\rho_\star}
 +  \tau \int_\Omega   H_\eta(\linstress \fp \linstrain(\bv) )H_\epsilon (|\linstress_\delta|^2 - \kappa_\star^2) \linstrain(\bv) \fp \linstress \\
&\geq 
 \|\linstress\|^2_{\cA}
+   \|\bv\|^2_{\rho_\star}
- \frac{1}{2}\|\linstress^{k-1}\|^2_{\cA}
- \frac{1}{2} \|\linstress\|^2_{\cA}
- \|\bv^{k-1}\|^2_{\rho_\star}
- \frac{1}{2}  \|\bv\|^2_{\rho_\star}
- \frac{\tau}{2}  \|\bv\|^2_{\rho_\star}
- \frac{\tau}{2} \|\mathbf{f}^k\|^2_{\rho_\star} \\
&\quad 
- \frac{c_{\textrm{tr}}^2 \tau}{2\rho_\star^- h}\|\bv\|^2_{\rho_\star}
- \frac{\tau}{2}\|\mathbf{t}^k_b\|_{L^2(\partial\Omega_N)}
- \frac{c_{\textrm{inv}}\tau}{h (A_{\textrm{min}} \rho^-_\star)^{1/2}}\|\bv\|_{\rho_{\star}} \|\linstress\|_{\cA} \\
&\geq 
\frac{1}{2}\left(1 - \frac{c_{\textrm{inv}}}{(A_{\textrm{min}} \rho^-_{\star})^{1/2} }\frac{\tau}{h}\right)
 \|\linstress\|^2_{\cA}
+ 
\frac{1}{2}\left(1 - \frac{ (A_{\textrm{min}} \rho^-_{\star})^{1/2} + A_{\textrm{min}}^{1/2}c^2_{\textrm{tr}} +  c_{\textrm{inv}}{\rho^-_\star}^{1/2} }{A_{\textrm{min}}^{1/2} \rho^-_{\star} }\frac{\tau}{h}\right)
\|\bv\|^2_{\rho_\star}\\
&\quad 
- \frac{1}{2}\|\linstress^{k-1}\|^2_{\cA}
- \|\bv^{k-1}\|^2_{\rho_\star}
- \frac{1}{2} \|\mathbf{f}\|^2_{L^2_{\rho_\star}(Q)}
+ \frac{1}{2}\|\mathbf{t}_b\|_{L^2(0,T;L^2(\partial\Omega_N))},
\end{align*}
\endgroup
where we employed Young's inequality, the inverse inequalites \eqref{eq:inverse_inequalities}, and the fact that $\tau \|\mathbf{f}^k\|^2_{\rho_\star} \leq \|\mathbf{f}\|^2_{L^2_{\rho_\star}(\Omega)}$. Hence, the claim follows if we assume that 
\begin{equation}\label{eq:smallness}
\frac{\tau}{h} < \frac{A_{\textrm{min}}^{1/2}\rho^-_\star}{c_{\textrm{inv}}{\rho^-_\star}^{1/2} + c^2_{\textrm{tr}} A_{\textrm{min}}^{1/2} + A^{1/2}_{\textrm{min}} {\rho^-_{\star}}^{1/2}}.
\end{equation}
The same corollary to Brouwer's fixed point theorem in addition implies that the solution is bounded,
\begin{equation*}
\|\linstress^k_{\eta}\|^2_{\cA} + \|\bv^k_\eta\|^2_{\rho_\star} \leq \hat{c}.
\end{equation*}
We remark that it is likely that existence of discrete solutions can be proved without assuming a condition like \eqref{eq:smallness} by relying on the equivalence of norms in finite dimensions and the fact that $\tau$ and $h$ are fixed. However, we choose to stick to the argument presented above, since the condition \eqref{eq:smallness} will appear once again in the next section where we analyse the stability of the numerical scheme, for which uniform bounds are desirable.

Now, since the bounds are independent of the regularisation parameter $\eta$ in the first Heaviside function $H_\eta$, the Heine--Borel theorem implies that up to a subsequence, for every $k\in \{1,\ldots,T/\tau\}$ the sequence of solutions $ \linstress^k_\eta$ (here we make the $\eta$-dependence explicit) converges as $\eta\to 0$,
\begin{align*}
  \linstress^k_\eta &\to \linstress^k &&\textrm{strongly in }L^\infty(\Omega)^{d\times d}, \\
  \bv^k_\eta &\to \bv^k &&\textrm{strongly in }W^{1,\infty}(\Omega)^{d},
\end{align*}
for some $\linstress^k\in \Sigma^h$ and $\bv^k\in V^h$. At this point we have used the fact that weak and strong convergence are equivalent in finite-dimensional spaces. This implies in particular that $H_\eta(\linstress^k_\eta \fp \linstrain(\bv^k_\eta)) \to H(\linstress^k\fp\linstrain(\bv^k))$ pointwise a.e.\ in $\Omega$, and so the limiting functions satisfy the system with the unregularised Heaviside function.

In summary, numerical solutions are guaranteed to exist, assuming the ratio $\frac{\tau}{h}$ is small enough. We note also that very similar arguments yield existence of solutions for the displacement-temperature system \eqref{eq:discrete_ut}.

\begin{remark}
We could also consider semi-implicit discretisation schemes such as
\begin{equation}
\begin{aligned}
\int_\Omega \cA(\dtau \linstress^k)\fp \bbsigma
- \int_\Omega \linstrain(\bv^k)\fp\bbsigma
+\int_\Omega &  H(\linstress^k \fp \linstrain(\bv^k) )H_\epsilon (|\linstress_\delta^k|^2 - \kappa_\star^2) \linstrain(\bv^k) \fp \bbsigma = 0   &&\forall \bbsigma \in \Sigma^h,\\
\int_\Omega\rho_\star \dtau \bv^k \cdot \bw 
               +\int_\Omega \linstress^{k-1} &\fp \linstrain(\bw) = \int_\Omega \rho_\star \mathbf{f}^{k-1}\cdot \bw
+ \int_{\partial\Omega_N} \mathbf{t}^{k-1} \cdot \bw
                             && \forall \bw\in V^h,
\end{aligned}
\end{equation}
and a similar analysis applies. The difference in this scheme compared to \eqref{eq:discrete_Tv} is that here the velocity $\bv^k$ is computed first using the information at time $t^{k-1}$ and the equation for $\linstress^k$ is solved afterwards. In the absence of plastic behaviour this results in a symplectic scheme that conserves a (modified) energy.
\end{remark}

\begin{remark}
A consequence of the fact that $\BD(V^h) \subset \Sigma^h$ and that $V^h\subset H^1_{\partial\Omega _D}(\Omega)^d$ is that the following discrete inf-sup condition holds:
\begin{equation}\label{eq:infsup}
\adjustlimits 
\inf_{\bw\in V^h} \sup_{\bbsigma\in \Sigma^h}\frac{\int_\Omega \bbsigma\fp \linstrain(\bw)}{\|\bw\|_{H^1(\Omega)}\|\bbsigma\|^2_{L^2(\Omega)}}
\geq 
\gamma_\star>0.
\end{equation}
where $\gamma_\star>0$ is independent of $h$; note that the above is then simply a reformulation of Korn's inequality. The validity of \eqref{eq:infsup} is not essential for the analysis of the discrete problem \eqref{eq:discrete_Tv}, but it would be crucial if we were interested in solving the quasi-static problem (i.e.\ without the time derivative $\dotd{\bv}$).
\end{remark}

\subsection{Stability}
Now we take a more careful look at the stability of the scheme. Let us first look at the continuous system \eqref{eq:pde_Tv}. First we assume that solutions are smooth enough so that all subsequent manipulations are well-defined in the classical sense. The  multiplication of the first equation in system~\eqref{eq:pde_Tv} by $\bv$, the second by $\linstress$ and integrating over $\Omega$ results in the energy balance in the form
\begin{equation}
  \frac{1}{2}\frac{\mathrm{d}}{\mathrm{d}t}\left(\int_\Omega \cA(\linstress)\fp \linstress + \rho_\star|\bv|^2 \right)
  + \int_\Omega H(\linstress \fp \linstrain(\bv) )H_\epsilon (|\linstress_\delta|^2 - \kappa_\star^2)\linstress\fp \linstrain(\bv)
  =
  \int_\Omega \rho_\star \mathbf{f}\cdot \bv
  +\int_{\partial\Omega_N} \mathbf{t}_b\cdot \bv.
\end{equation}
If we follow the terminology used in~\cite{SH.2006}, thus, if we denote the kinetic energy by $E_{\textrm{kin}}(\bv):= \int_\Omega\frac{1}{2}\rho_\star |\bv|^2$, the elastic potential energy by $E_{\textrm{int}}(\linstress):=\int_\Omega \frac{1}{2}\cA(\linstress)\fp \linstress$, and the potential energy associated with the applied loads by $E_{\textrm{ext}}(\bu):= -\int_\Omega \rho_\star \mathbf{f}\cdot \bu - \int_{\partial\Omega_N}\mathbf{t}_b\cdot \bu$, then the energy balance can be rewritten as
\begin{equation}\label{eq:energy_balance2}
  \frac{\mathrm{d}}{\mathrm{d}t}\left[ E_{\textrm{kin}}(\bv) + E_{\textrm{int}}(\linstress) + E_{\textrm{ext}}(\bu) \right]
  = -\int_\Omega H(\linstress \fp \linstrain(\bv) )H_\epsilon (|\linstress_\delta|^2 - \kappa_\star^2)\linstress\fp \linstrain(\bv)
  \leq 0.
\end{equation}
Inspecting~\eqref{eq:energy_balance2}, it is clear that there is mechanical dissipation only when the material is being loaded and the yield stress has been reached. Moreover, if we define the thermal energy as $E_{\textrm{th}}(\theta) := \int_\Omega \rho_\star c_v \theta$, then integrating the temperature equation \eqref{eq:pde_temperature} and adding the result to \eqref{eq:energy_balance2}, yields the total energy balance

\begin{equation}\label{eq:energy_balance_total}
\frac{\mathrm{d}}{\mathrm{d}t}\left[ E_{\textrm{kin}}(\bv) + E_{\textrm{int}}(\linstress) + E_{\textrm{ext}}(\bu) + E_{\textrm{th}}(\theta) \right]=0.
\end{equation}
The balance \eqref{eq:energy_balance_total} highlights the fact that, as a consequence of the thermodynamically consistent derivation of the model, all the various energy dissipation mechanisms are accounted for in the model.

We now obtain an analogue of \eqref{eq:energy_balance2} at the \emph{discrete level}. Choosing $\bbsigma:=\linstress^k$ and $\bw=\bv^k$ in the numerical formulation \eqref{eq:discrete_Tv}, and using the elementary identity $(a-b)a = \frac{1}{2}a^2 - \frac{1}{2}b^2 + \frac{1}{2}|a-b|^2$ for two numbers $a,b\in \R$, yields for all $k\in \{1,\ldots,T/\tau\}$ the equality
\begin{align*}
&\frac{1}{2\tau} \|\bv^k\|^2_{\rho_\star}
-\frac{1}{2\tau} \|\bv^{k-1}\|^2_{\rho_\star}
+ \frac{1}{2\tau} \|\linstress^k\|^2_{\cA}
- \frac{1}{2\tau} \|\linstress^{k-1}\|^2_{\cA}
+ \frac{1}{2\tau} \|\bv^k-\bv^{k-1}\|^2_{\rho_\star} \\
& + \frac{1}{2\tau} \|\linstress^k - \linstress^{k-1}\|^2_{\cA}
+  \int_\Omega   H(\linstress^k \fp \linstrain(\bv^k) )H_\epsilon (|\linstress^k_\delta|^2 - \kappa_\star^2) \linstrain(\bv^k) \fp \linstress^k
= \int_\Omega \rho_\star \mathbf{f}\cdot \bv^k 
+ \int_{\partial\Omega_N}\mathbf{t}_b\cdot \bv^k.
\end{align*}
Hence, if we define the numerical dissipation $\mathcal{D}_{\tau}^k := \frac{1}{2\tau}\|\bv^k - \bv^{k-1}\|^2_{\rho_\star} + \frac{1}{2\tau}\|\linstress^k - \linstress^{k-1}\|^2_{\cA}$, then the equality just derived above can be rewritten as 
\begin{equation}\label{eq:energy_balance_discrete_a}
  \dtau\left[E_{\textrm{kin}}(\bv^k) + E_{\textrm{int}}(\linstress^k) \right]
  + E_{\textrm{ext}}(\bv^k) 
  = - \mathcal{D}_\tau^k
-  \int_\Omega   H(\linstress^k \fp \linstrain(\bv^k) )H_\epsilon (|\linstress^k_\delta|^2 - \kappa_\star^2) \linstrain(\bv^k) \fp \linstress^k
\leq 0.
\end{equation}
This equality clearly mimics the continuous energy balance \eqref{eq:energy_balance2}, except for the presence of the numerical dissipation term $\mathcal{D}^k_\tau$. Moreover, testing the temperature equation \eqref{eq:discrete_ut} with $\phi=1$ we also obtain the discrete total energy balance:
\begin{equation}\label{eq:energy_balance_discrete}
  \dtau\left[E_{\textrm{kin}}(\bv^k) + E_{\textrm{int}}(\linstress^k) + E_{\textrm{th}}(\theta^k) \right]
  + E_{\textrm{ext}}(\bv^k) 
  = - \mathcal{D}_\tau^k
\leq 0,
\end{equation}
which is analogous to the total energy balance \eqref{eq:energy_balance_total}, up to numerical dissipation.

\begin{remark}
If we denote the piecewise-linear (in time) interpolant of the sequence $\{\bv^k\}_{k=0}^{T/\tau}$ by $\tilde{\bv}_{h,\tau} \in C([0,T];V^h)$, then we see that
\begin{equation*}
\frac{1}{\tau}\|\bv^k-\bv^{k-1}\|^2_{L^2(\Omega)} 
 = \tau \left\|\frac{\partial \tilde{\bv}_{h,\tau}(t^k_-)}{\partial t}\right\|^2_{L^2(\Omega)},
\end{equation*}
and so if the problem satisfies appropriate regularity properties so that the norm on the right-hand-side is bounded, then it is clear that the numerical dissipation term vanishes as $\tau\to 0$. Similar arguments apply to the stress.
\end{remark}

Now, let us denote the piecewise-constant (in time) interpolant associated to the sequence $\{\bv^k\}_{k=0}^{T/\tau}$ by $\bv_{h,\tau} \in L^\infty((0,T);V^h)$, and define $\linstress_{h,\tau}$ analogously. Then, multiplying the discrete energy balance \eqref{eq:energy_balance_discrete} by $2\tau$, using a similar argument to the one employed in the previous section, and summing over $k$, we obtain the stability estimate

\begin{equation}\label{eq:discrete_stability}
\|\bv_{h,\tau}\|^2_{L^\infty(0,T;L_{\rho_\star}^2(\Omega))}
+ \|\linstress_{h,\tau}\|^2_{L^\infty(0,T;L^2_{\cA}(\Omega))}
+ \sum_{k=1}^{T/\tau} \tau \mathcal{D}^k_\tau 
\leq 
\|\bv_0\|^2_{\rho_\star}
+ \|\linstress_0\|^2_{\cA}
+ \|\mathbf{f}\|^2_{L_{\rho_\star}^2(Q)}
+ \|\mathbf{t}_b\|^2_{L^2(0,T;L^2(\partial\Omega_N))},
\end{equation}
where we assume that the condition \eqref{eq:smallness} is satisfied. We remark here that analogous arguments apply to the displacement and the temperature system~\eqref{eq:discrete_ut}.

\begin{remark}\label{rmk:infsup_bound}
From the inf-sup condition \eqref{eq:infsup} we can also try to obtain a bound for the discrete velocities
\begin{equation*}
\tau\|\bv^k\|_{H^1(\Omega)} \leq
\tau \sup_{\bbsigma\in \Sigma^h}\frac{\int_\Omega \linstrain(\bv^k)\fp \bbsigma}{\|\bbsigma\|_{L^2(\Omega)}}
\leq
A_{\mathrm{min}}^{-1}(\|\linstress^k\|_{\cA} + \|\linstress^{k-1}\|_{\cA}) 
+ \tau \|\linstrain(\bv^k)\|_{L^2(\Omega)}.
\end{equation*}
In the absence of plastic behaviour the last term is not present and this would imply, together with \eqref{eq:discrete_stability}, that we can bound uniformly $\|\bv_{h,\tau}\|_{L^2(0,T;H^1(\Omega))}$ in terms of the data; this is what would be expected in the linear elasticity model. However, in general this only yields a bound for $\bv_{h,\tau}$ in $L^2(Q)^d$, which does not improve \eqref{eq:discrete_stability}. This lack of a priori boundedness of the velocity gradients is what results in the restriction on the ratio $\frac{\tau}{h}$.
\end{remark}

\begin{remark}\label{rmk:convergence}
Note that
\begin{equation*}
 \sum_{k=1}^{T/\tau} \tau \mathcal{D}^k_\tau 
= \sum_{k=1}^{T/\tau} \|\bv^k - \bv^{k-1}\|^2_{\rho_\star}
+ \|\linstress^{k}- \linstress^{k-1}\|^2_{\cA}
= \|\partial_t \bv_{h,\tau}\|_{\mathcal{M}(0,T;L_{\rho_\star}^2(\Omega))}
+ \|\partial_t \linstress_{h,\tau}\|_{\mathcal{M}(0,T;L_{\cA}^2(\Omega))},
\end{equation*}
and so the discrete stability estimate \eqref{eq:discrete_stability} also yields a bound on the time derivatives of the approximate solutions; here $\mathcal{M}(0,T;L_{\rho_\star}^2(\Omega))$ denotes the space of Radon measures in time with values into $L^2_{\rho_\star}(\Omega)$. (The space $\mathcal{M}(0,T;L_{\cA}^2(\Omega))$ is defined analogously.) This is enough, for example by applying \cite[Cor.\ 7.9]{Rou.2013}), to prove that as $\tau \to 0$, the solutions $(\linstress_{h,\tau},\bv_{h,\tau})$ converge to functions $(\linstress_h,\bv_h)$ that solve the system
\begin{equation}
  \begin{aligned}
    \rho_\star \dotd{\bv_h} &-  \diver\linstress_h 
    = \rho_\star\mathbf{f} \qquad \quad &&\text{ in }{V^h},\\
    \cA(\dotd{\linstress}_h) = \linstrain(\bv_h)
                          &- H(\linstress_h\fp \linstrain(\bv_h) )H_\epsilon (|(\linstress_h)_\delta|^2 - \kappa_\star^2)\linstrain(\bv_h)
                          &&\text{ in }{\Sigma^h}.
  \end{aligned}
\end{equation}
Obtaining convergence as $h\to 0$ is a more delicate matter given the relatively weak bounds available to us (c.f.\ Remark \ref{rmk:infsup_bound}). In fact, at this point we face the lack of analytical results for a system of partial differential equations of the rate-type~\eqref{eq:pde_Tv}---it is not completely obvious what the proper notion of weak solution should be. Conceivably, this problem could become more tractable by introducing hardening into the model, and then the solutions for the perfect plasticity model would be obtained in a vanishing hardening limit; this will be the subject of future research.
\end{remark}

{\color{black}
\begin{remark}\label{rmk:temp_dependent_parameters}
The proofs of existence of discrete solutions and discrete stability rely mainly on two ingredients: continuity of the residual~\eqref{eq:discrete_residual} and \emph{a priori} estimates. As a consequence, incorporating temperature-dependent material parameters into the analysis is straightforward as long as these properties are not affected. For instance, the results apply also to the model with a temperature-dependent yield stress $\kappa_\star(\theta)$, assuming that this dependence is continuous (and since it appears exclusively in the argument of a Heaviside function, boundedness of $\kappa_\star$ is not required); temperature-dependence in the compliance tensor $\mathcal{A}$ is also fine, as long as the norm $\|\cdot\|_{\cA}$ remains equivalent to the $L^2$-norm, which is the case whenever the material parameters $\ymod=\ymod(\theta)$ and $\nu=\nu(\theta)$ are continuous functions of the temperature, which are bounded from above and below by positive constants. (This, together with a more complex temperature evolution equation, would constitute the first step in analysis of more complex models that take into account for example the Gough--Joule effect.) Evidently, with such a temperature-dependence, the mechanical and thermal systems \eqref{eq:discrete_Tv} and \eqref{eq:discrete_ut} do not decouple and have to be solved simultaneously.
\end{remark}
}

\section{Numerical experiments}
We now implement the discrete formulations \eqref{eq:discrete_Tv} and \eqref{eq:discrete_ut} to illustrate that they indeed capture the behaviour expected from the model. (Note that we are not using a quasistatic approximation, both equations are treated as evolutionary equations. In particular the acceleration term in the balance of linear momentum is taken into account in the proposed numerical algorithm as well as in its implementation.) We first implement the problem in one spatial dimension, and we show that the mechanical response is as expected during one loading-unloading cycle. Subsequently we implement the problem describing a two dimensional plate with an elliptical hole. The nonlinear systems for the stress and velocity at each time step are handled with Newton's method supplemented with the error oriented line search \texttt{NLEQERR} from PETSc \cite{PETSc}; the absolute and relative tolerances for the nonlinear solver are set to $10^{-6}$. The linear systems at each Newton step are solved using the LU factorisation algorithm from MUMPS \cite{MUMPS:1}; the linear systems for the displacement and the temperature are solved in turn using MUMPS as well. Everything is implemented through the finite element software \texttt{firedrake} \cite{Firedrake}; the code used to implement the computational experiments, including the exact components of \texttt{firedrake} that have been employed, has been archived in Zenodo (\url{https://zenodo.org/record/7342357}) \cite{zenodo} for reproducibility purposes.

\subsection{One dimensional mechanical response}
We solve the problem on the unit interval $\Omega = (0,1)$ and for times $t\in [0,1]$. (If not stated otherwise all physical quantities are given in the SI base units or use combination thereof.) We impose boundary conditions on the displacement:
\begin{equation}
u(t,0) = -u(t,1) :=
\left\{
  \begin{array}{cc}
    -\frac{1}{10} e^{1 + \frac{1}{4t(t-1)}} & t\in (0,1).\\
    0 & \text{otherwise}.
\end{array}
\right.
\end{equation}
This describes loading for $t\in (0, \tfrac{1}{2})$\,s and unloading otherwise. Since the problem is one-dimensional, we denote the scalar displacement, stress and strain are denoted by $u$, $\sigma$, and $\varepsilon$, respectively. We set the Young modulus to $\ymod = 10^4$\,Pa.
(The values of physical constants are in this example entirely artificial.) We employ a simple continuation algorithm with respect to $\epsilon$ to produce better initial guesses for Newton's method; for instance, the problem is solved with a larger (and thus easier) value for $\epsilon$ and the solution is used as an initial guess for the problem with regularisation parameter $\epsilon - \delta \epsilon$ until the desired value is reached.

\begin{figure}
\centering
\subfloat[C][{\centering Stress-strain response}]{{%
	\includegraphics[width=0.5\textwidth]{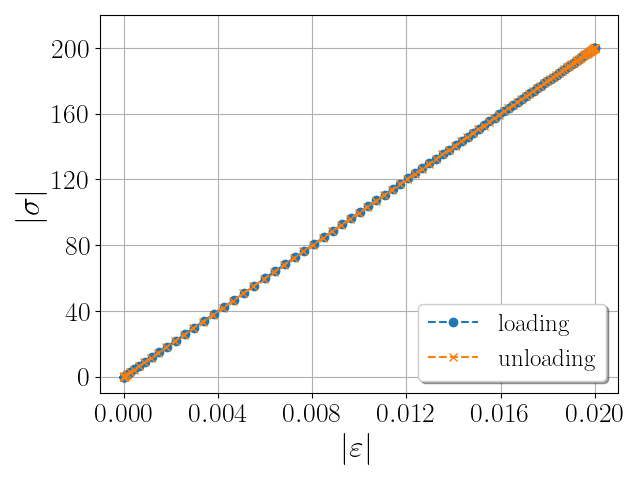}%
	}}%
	\subfloat[D][{\centering Maximum stress}]{{%
	\includegraphics[width=0.5\textwidth]{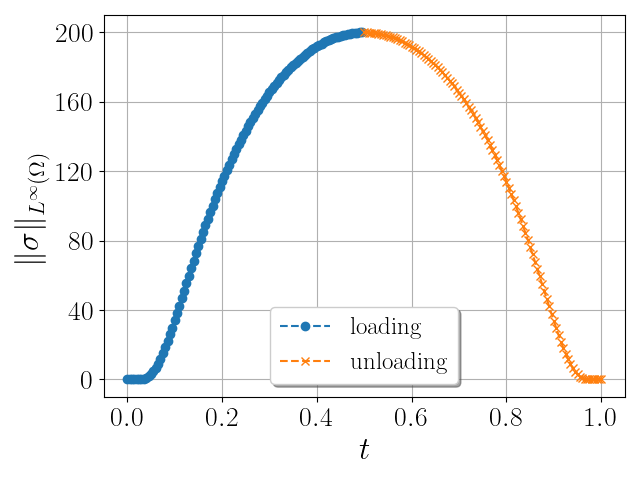}%
	}}\\
	\caption{Mechanical response at $X=0.75$\,m for the problem with $\kappa_\star = 10^7$\,Pa.}%
	\label{fig:1D_elastic}
\end{figure}

Figure \ref{fig:1D_elastic} shows the stress-strain response at the point $X=0.75$\,m with a very large yield-stress $\kappa_\star = 10^7$\,Pa; the problem is solved with 482 spatial degrees of freedom and a time step $\tau = 5\times 10^{-4}$\,s; a plot of the maximum stress $\|\sigma\|_{L^\infty(\Omega)}$ with respect to time is also shown for reference. This is simply a sanity check to verify that the solution of our proposed numerical scheme behaves as expected; namely, the solution exhibits solely elastic behaviour.

\begin{figure}
\centering
\subfloat[C][{\centering Heaviside function $H_\epsilon^{(1)}$, $\epsilon = 10$}]{{%
	\includegraphics[width=0.5\textwidth]{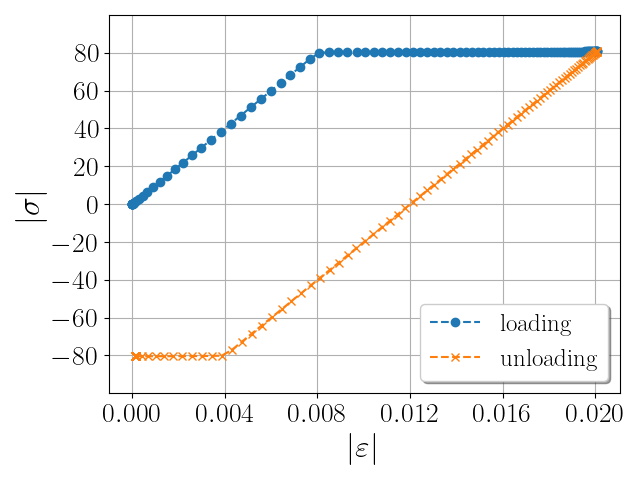}%
	}}%
	\subfloat[D][{\centering Heaviside function $H_\epsilon^{(1)}$, $\epsilon = 200$}]{{%
	\includegraphics[width=0.5\textwidth]{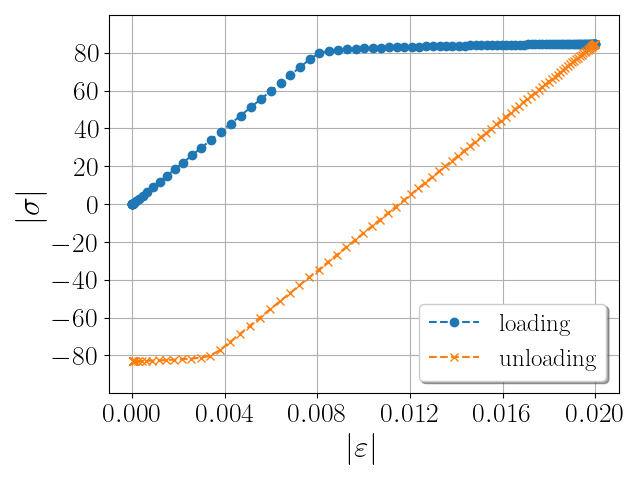}%
	}}\\
	\caption{Stress-strain response at $X=0.75$\,m for the problem with $\kappa_\star = 80$\,Pa.}%
	\label{fig:1D_plastic_stress-strain}
\end{figure}

\begin{figure}
\centering
\subfloat[C][{\centering Heaviside function $H_\epsilon^{(1)}$, $\epsilon = 10$}]{{%
	\includegraphics[width=0.5\textwidth]{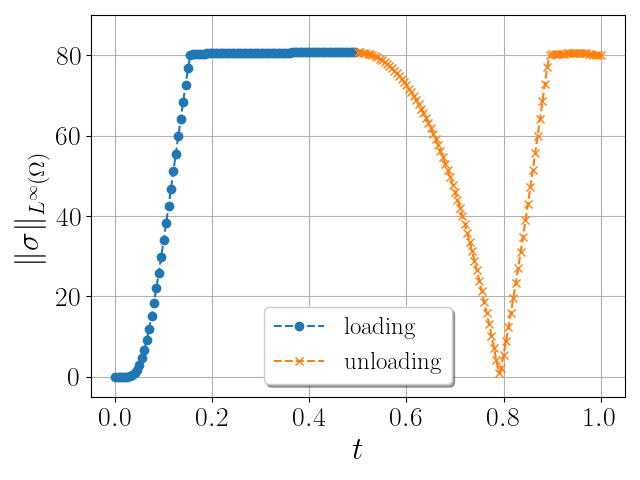}%
	}}%
	\subfloat[D][{\centering Heaviside function $H_\epsilon^{(1)}$, $\epsilon = 100$}]{{%
	\includegraphics[width=0.5\textwidth]{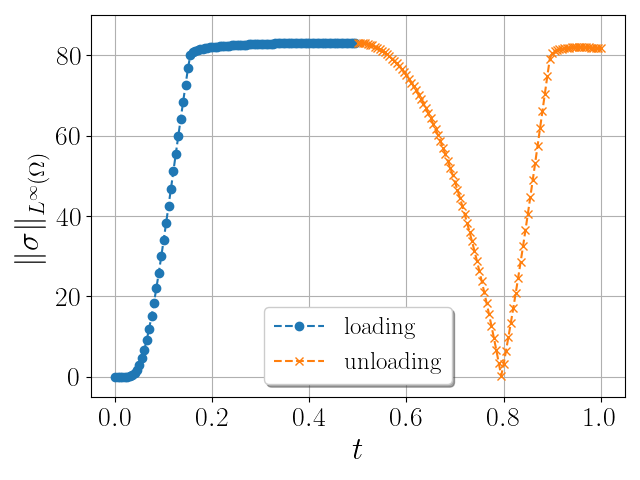}%
	}}\\
\subfloat[C][{\centering Heaviside function $H^{(2)}_\epsilon$, $\epsilon = 10$}]{{%
	\includegraphics[width=0.5\textwidth]{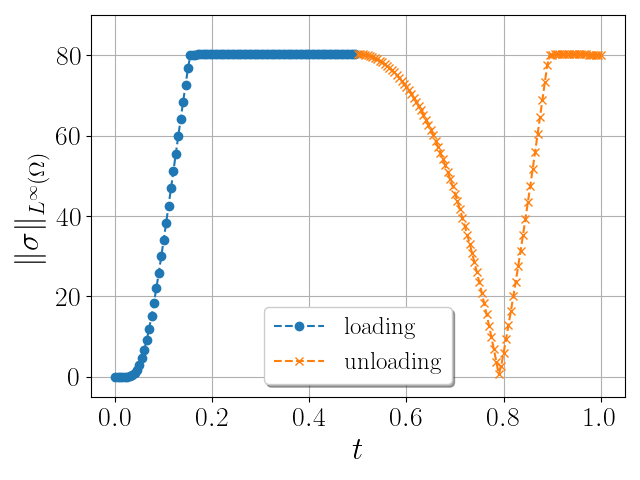}%
	}}%
	\subfloat[D][{\centering Heaviside function $H^{(2)}_\epsilon$, $\epsilon = 100$}]{{%
	\includegraphics[width=0.5\textwidth]{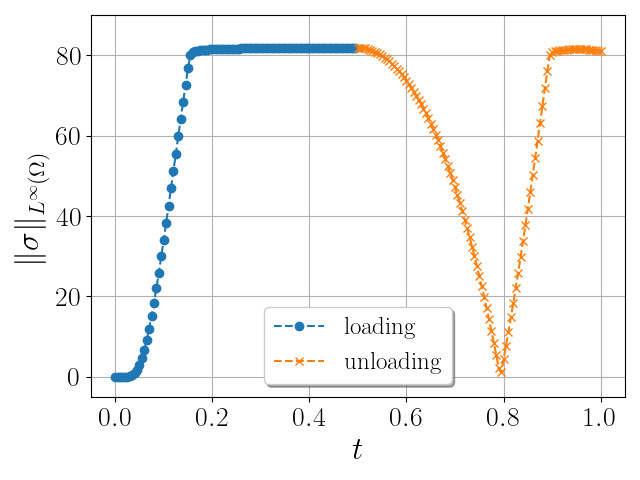}%
	}}\\
\subfloat[C][{\centering Heaviside function $H^{(3)}_\epsilon$, $\epsilon = 10$}]{{%
	\includegraphics[width=0.5\textwidth]{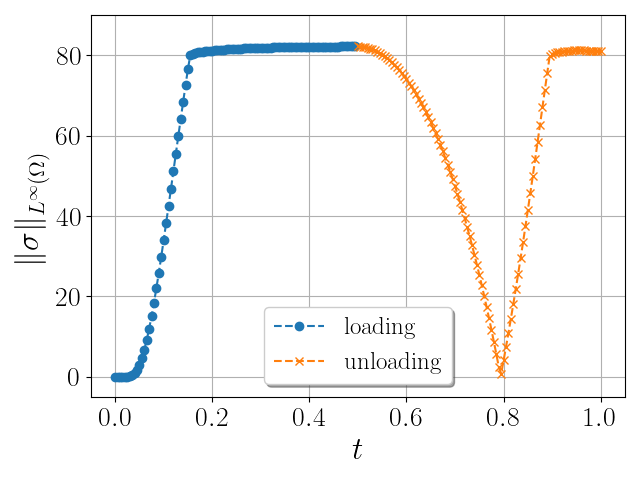}%
	}}%
	\subfloat[D][{\centering Heaviside function $H^{(3)}_\epsilon$, $\epsilon = 100$}]{{%
	\includegraphics[width=0.5\textwidth]{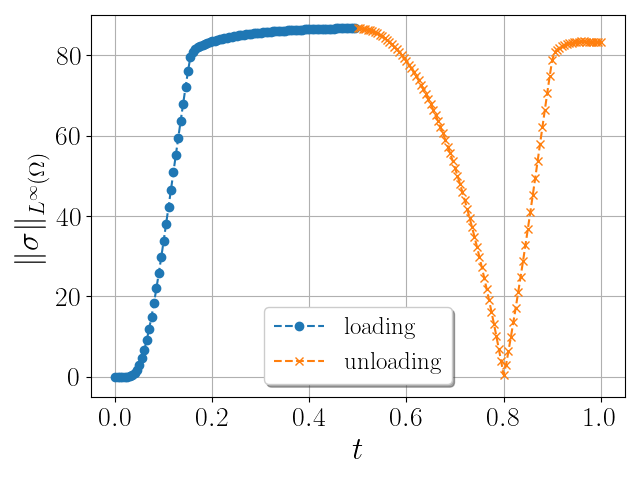}%
	}}\\
	\caption{Maximum stress over time for the problem with $\kappa_\star = 80$\,Pa, and various approximations of the Heaviside function.}%
	\label{fig:1D_plastic-maxstress}
\end{figure}

The same problem is subsequently solved for the yield stress of $\kappa_\star = 80$\,Pa with different values of $\epsilon$; the stress-strain relations are shown in Figure~\ref{fig:1D_plastic_stress-strain}. The values of the maximum stress are plotted in Figure \ref{fig:1D_plastic-maxstress} for different values of $\epsilon$ and the different approximations/regularisations of the Heaviside function. We observe that the numerical solutions capture the expected elastic--perfectly plastic behaviour  during one loading-unloading cycle. We also observe for large $\epsilon$ a non-sharp yield transition; depending on which Heaviside approximation we employ, the computed stress can be allowed to go slightly beyond $\kappa_\star$, but this effect disappears as $\epsilon$ decreases; in this regard we observe that the regularisation $H_\epsilon^{(2)}$ based on the hyperbolic tangent is the one that violates the constraint the least.

{
  \color{black}
As mentioned in Remark \ref{rmk:convergence}, a fully rigorous convergence analysis in this work is out of reach. However, we will supplement this example with convergence plots with respect to a reference solution, in order to illustrate what orders of convergence might be expected in practice; the reference solution was computed on a fine grid with 10240 elements in the spatial direction and 16384 elements in the temporal direction. From a given level of refinement level $l$, the level $l+1$ is obtained by reducing the (spatial and temporal) mesh size by half. The initial level $l=1$ consists of 40 and 64 elements in time and space, respectively.

As suggested from the stability estimate \eqref{eq:discrete_stability}, the natural way to measure the error is  the $L^\infty(0,T;L^2(\Omega))$. Namely, we compute the following errors for the stress and velocity, respectively:
\begin{equation}
E_{str}^l := \sup_{t\in [0,T]}\|\sigma^l - \sigma^*\|_{L^2(\Omega)}
\qquad
E_{vel}^l := \sup_{t\in [0,T]} \|v^l - v^*\|_{L^2(\Omega)},
\end{equation}
where $(\sigma^l,v^l)$ denotes the solution at level $l$ and $(\sigma^*,v^*)$ denotes the reference solution. Since the solutions are piecewise constants in time, the errors above can be computed in a straightforward manner. We note that to simplify the implementation we omit the weights corresponding to $\rho_\star$ and $\cA$ in the spatial norm; since these are positive, this does not affect the rate of convergence.

Figure \ref{fig:1D_convergence} shows the convergence behaviour for problems with various values for $\kappa_\star$ and $\epsilon$ (we employ the approximation $H^{(1)}_\epsilon$). For reference we plot as well the dependence on the number of spatial degrees of freedom $\mathrm{dofs}^\alpha$, with a certain convergence exponent $\alpha$; recall that the time step also decreases by half each level, so its influence is already correlated with the value of $\mathrm{dofs}$. We observe convergence with fractional rates, which is common for non-smooth problems with a possible loss of ellipticity (see e.g.\ \cite{WL.2011}).

\begin{figure}
\centering
\subfloat[C][{\centering Stress error, $\epsilon = 200$}]{{%
	\includegraphics[width=0.5\textwidth]{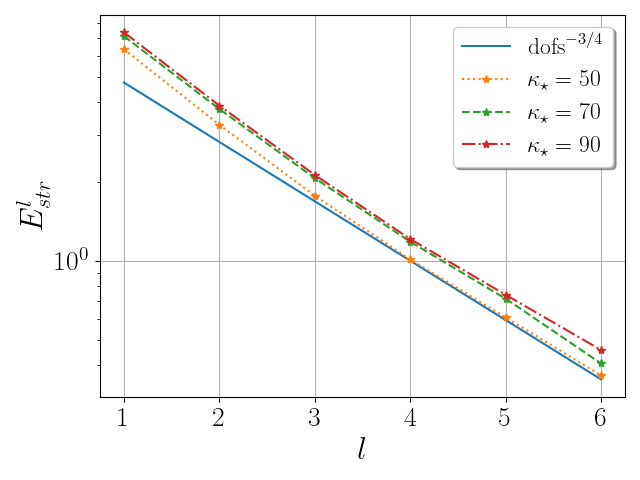}%
	}}%
	\subfloat[D][{\centering Velocity error, $\epsilon = 200$}]{{%
	\includegraphics[width=0.5\textwidth]{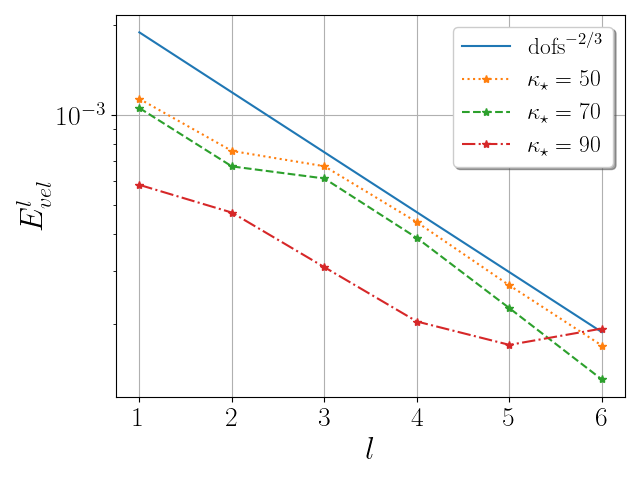}%
	}}\\
\subfloat[C][{\centering Stress error, $\kappa_\star = 70$\,Pa}]{{%
	\includegraphics[width=0.5\textwidth]{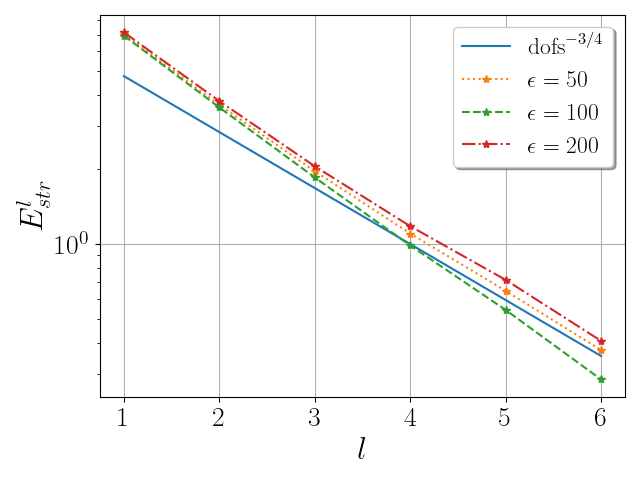}%
	}}%
	\subfloat[D][{\centering Velocity error, $\kappa_\star = 70$\,Pa}]{{%
	\includegraphics[width=0.5\textwidth]{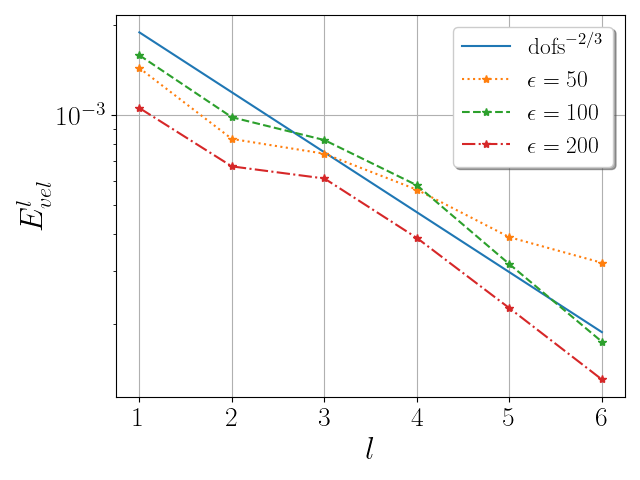}%
	}}\\
	\caption{Convergence behaviour for the solution of the 1D mechanical problem.}%
	\label{fig:1D_convergence}
\end{figure}
}

\subsection{2D plate with an elliptical hole}
\begin{figure}
\centering
	\includegraphics[width=0.45\textwidth]{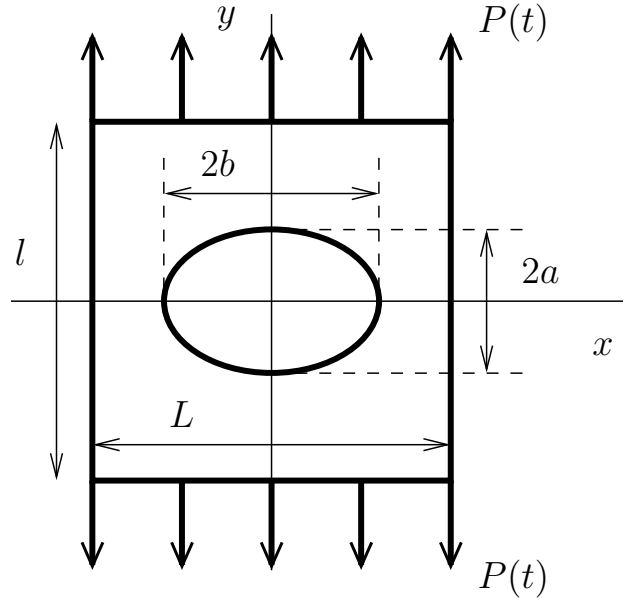}
	\caption{Square domain with an elliptical hole.}%
	\label{fig:ellipse}
\end{figure}

The problem is solved on the two-dimensional domain $(-\tfrac{L}{2},\tfrac{L}{2})\times (\tfrac{-l}{2},\tfrac{-l}{2})$, in which an elliptical hole is made, with semi-minor and semi-major axis of lengths $a$ and $b$, respectively (see Figure \ref{fig:ellipse}); in the implementation we choose $L=l=1$\,m, $a=0.3$\,m, and $b=0.5$\,m. We consider homogeneous Neumann boundary conditions for the temperature, while on the left and right boundaries and on the elliptical boundary we prescribe homogeneous natural boundary conditions for the mechanical problem, while on the top and bottom  we apply a time-dependent traction $P(t)$ in the vertical direction. This traction is meant to represent one loading-unloading cycle, and its magnitude taken to be a bump function
\begin{equation*}
P(t) := 
\left\{
  \begin{array}{cc}
    20 e^{1 + \frac{1}{4t(t-1)}} & t\in (0,1),\\
    0 & \text{otherwise}.
\end{array}
\right.
\end{equation*}

Concerning the initial conditions we assume that the body of interest is initially in the stress-free state and that the initial displacement is equal to zero. The initial temperature distribution is homogeneous in space, and in what follows we report only temperature changes with respect to this initial state.

Concerning the material parameters we take the values $\ymod = 10^4$\,Pa, $\nu=0.3$, $\kappa_\star = 1$\,Pa, $\kappa_{\textrm{th}}=1 \mathrm{W}\cdot \mathrm{m}^{-1}\cdot \mathrm{K}^{-1}$, and $c_v = 1\,\mathrm{J}\cdot \mathrm{kg}^{-1}\cdot \mathrm{K}^{-1}$.  (These parameter values are again artificial, we do not aim at a particular set of parameter values for a real material.) The problem is solved with a time step of $\tau = 5\times 10^{-4}$\,s, $1.52 \times 10^{5}$ degrees of freedom for the stress-velocity problem \eqref{eq:discrete_Tv} and $5.8 \times 10^{4}$ degrees of freedom for the displacement-temperature problem \eqref{eq:discrete_ut}. We employ the first Heaviside regularisation $H^{(1)}_\epsilon$ with the regularisation parameter fixed as $\epsilon=100$.

\begin{figure}
\centering
\subfloat[C][{\centering $t=0$\,s}]{{%
	\includegraphics[width=0.42\textwidth]{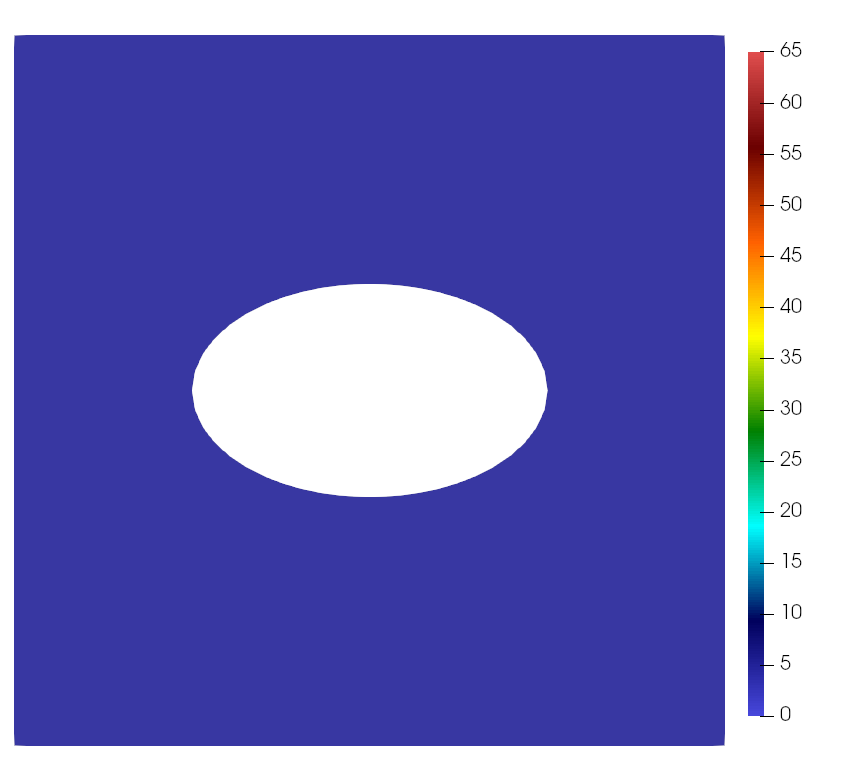}%
	}}%
	\subfloat[D][{\centering $t=0$\,s}]{{%
	\includegraphics[width=0.42\textwidth]{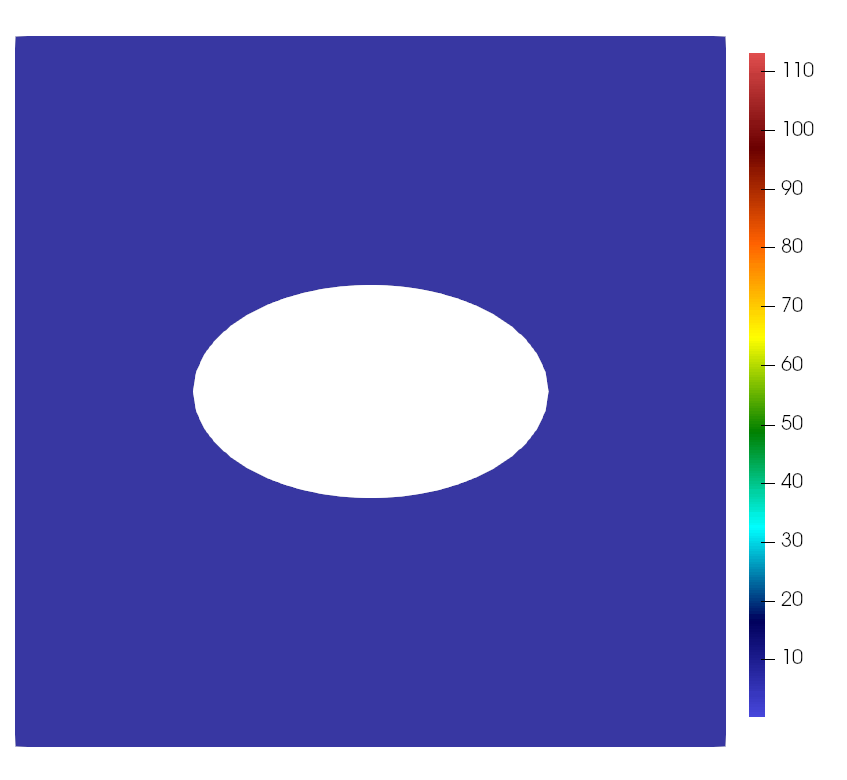}%
	}}\\
\subfloat[C][{\centering $t=0.25$\,s}]{{%
	\includegraphics[width=0.42\textwidth]{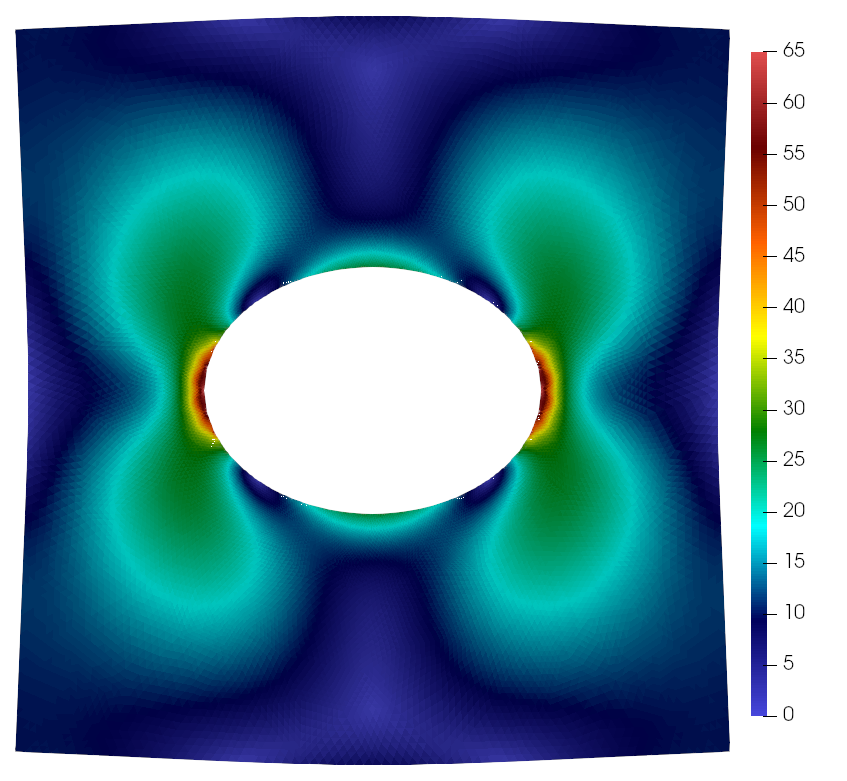}%
	}}%
	\subfloat[D][{\centering $t=0.25$\,s}]{{%
	\includegraphics[width=0.42\textwidth]{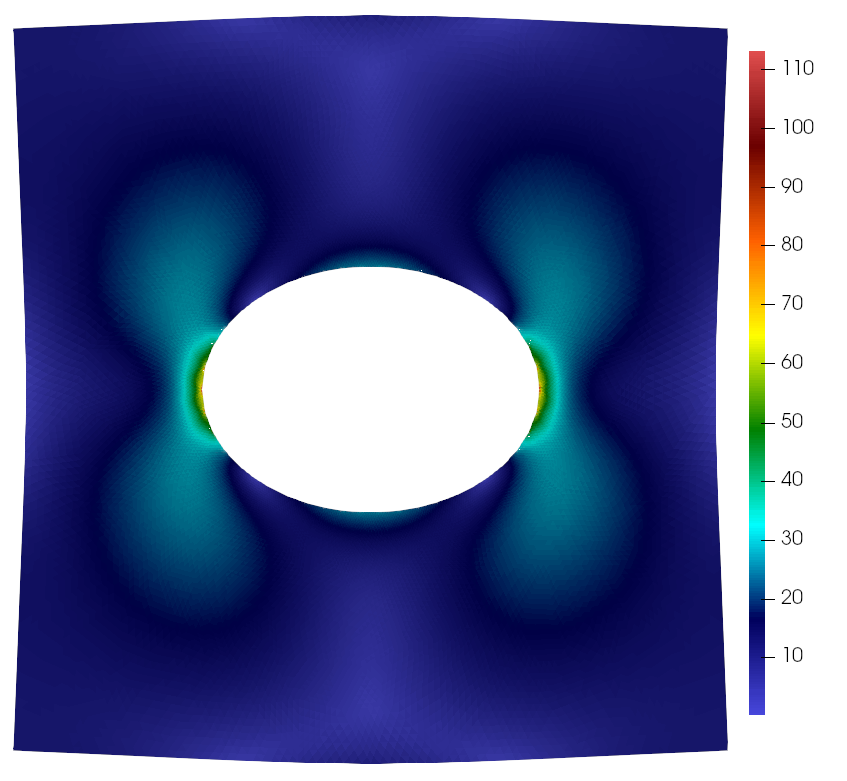}%
	}}\\
\subfloat[C][{\centering $t=0.5$\,s}]{{%
	\includegraphics[width=0.42\textwidth]{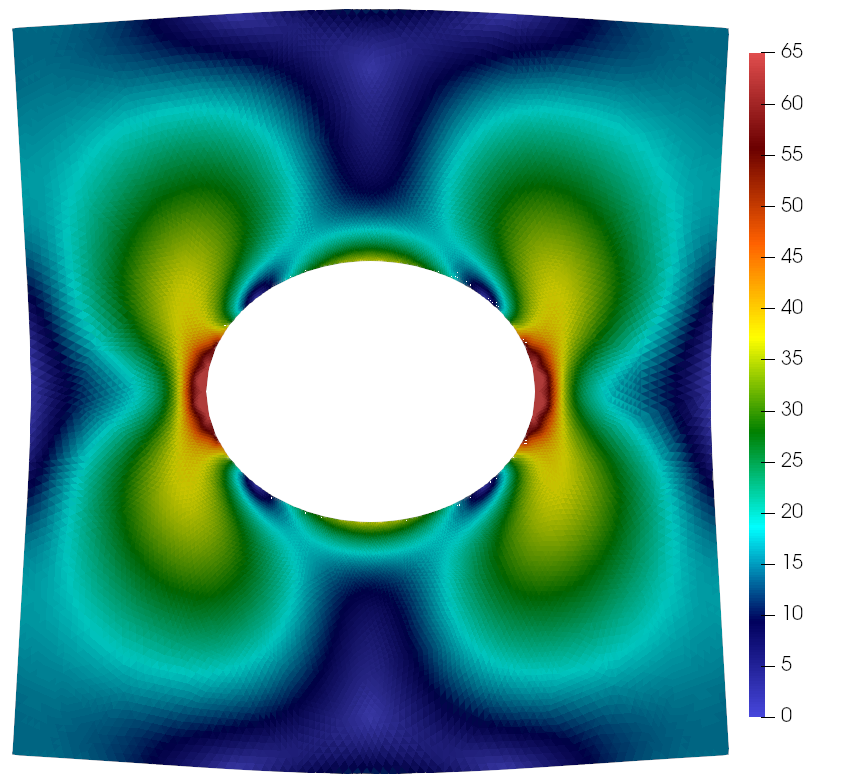}%
	}}%
	\subfloat[D][{\centering $t=0.5$\,s}]{{%
	\includegraphics[width=0.42\textwidth]{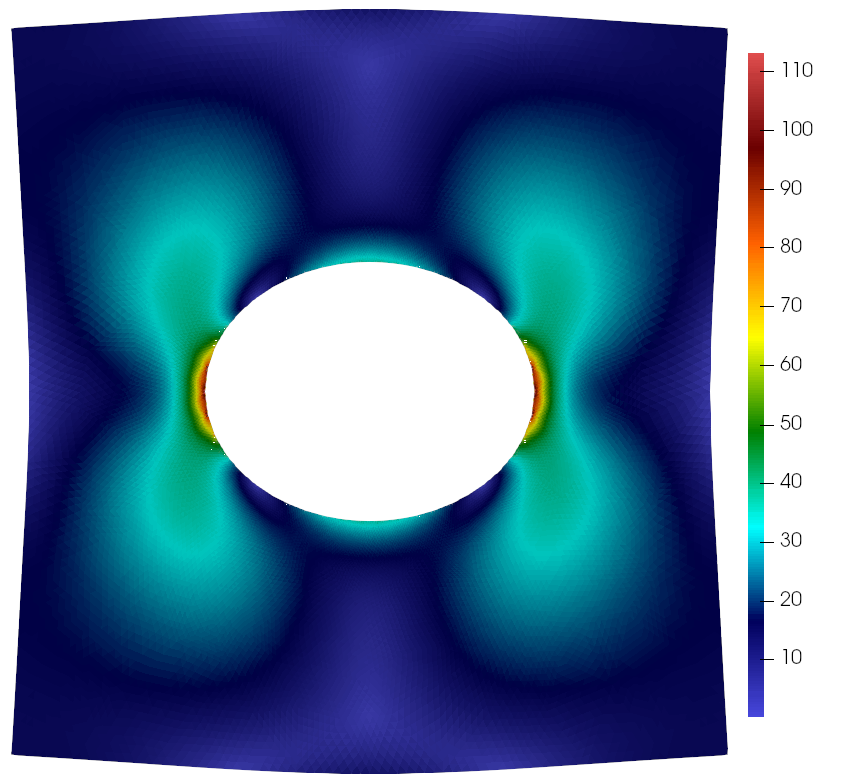}%
	}}\\
	\caption{Magnitude of $\linstress_\delta$ for the elastic problem with $\kappa_\star=10^7$\,Pa (right) and the problem with $\kappa_\star=60$\,Pa (left).}%
	\label{fig:ellipse_stress}
\end{figure}

\begin{figure}
\centering
\subfloat[C][{\centering $\kappa_\star = 60$\,Pa, $t=1$\,s}]{{%
	\includegraphics[width=0.45\textwidth]{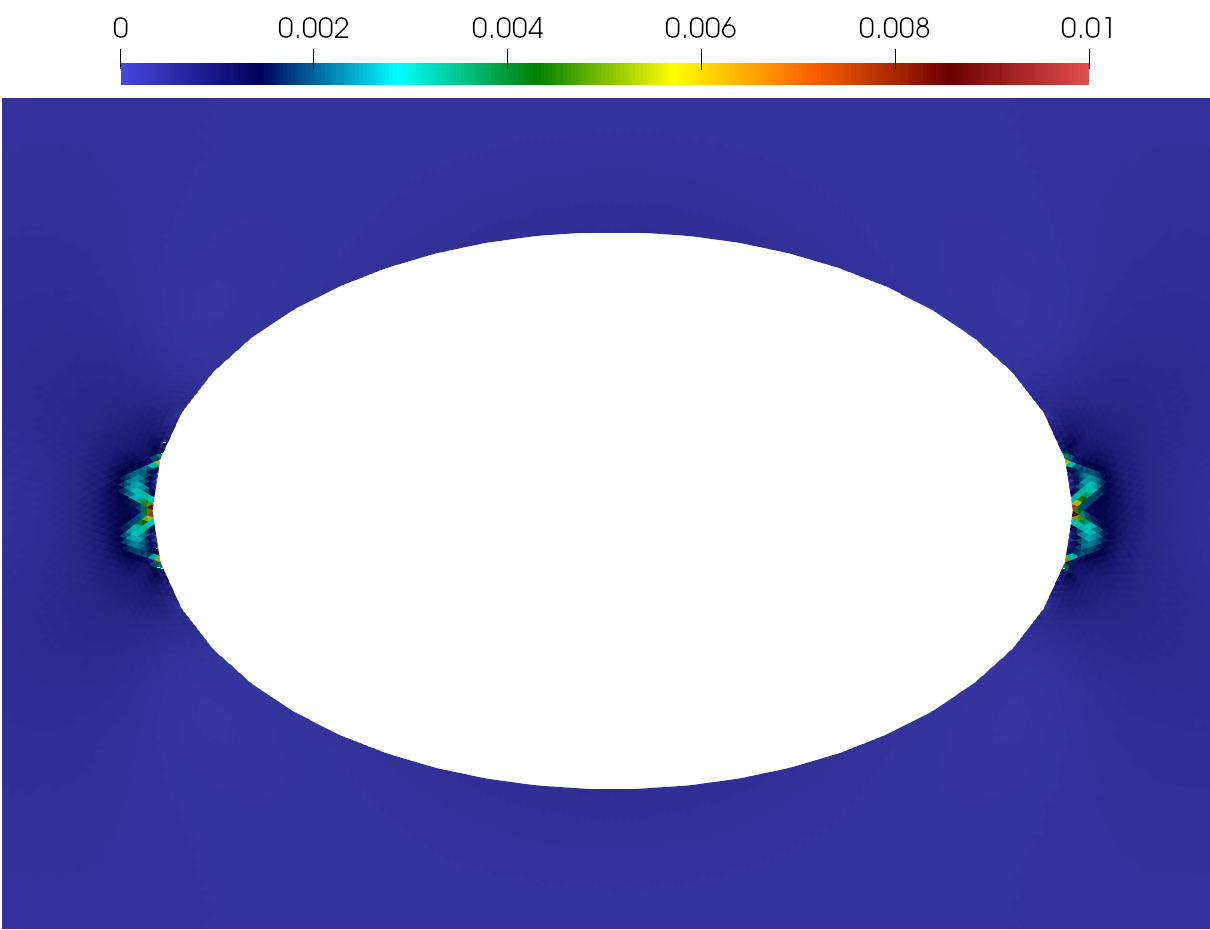}%
	}}%
\subfloat[D][{\centering $\kappa_\star = 10^7$\,Pa, $t=1$\,s}]{{%
	\includegraphics[width=0.45\textwidth]{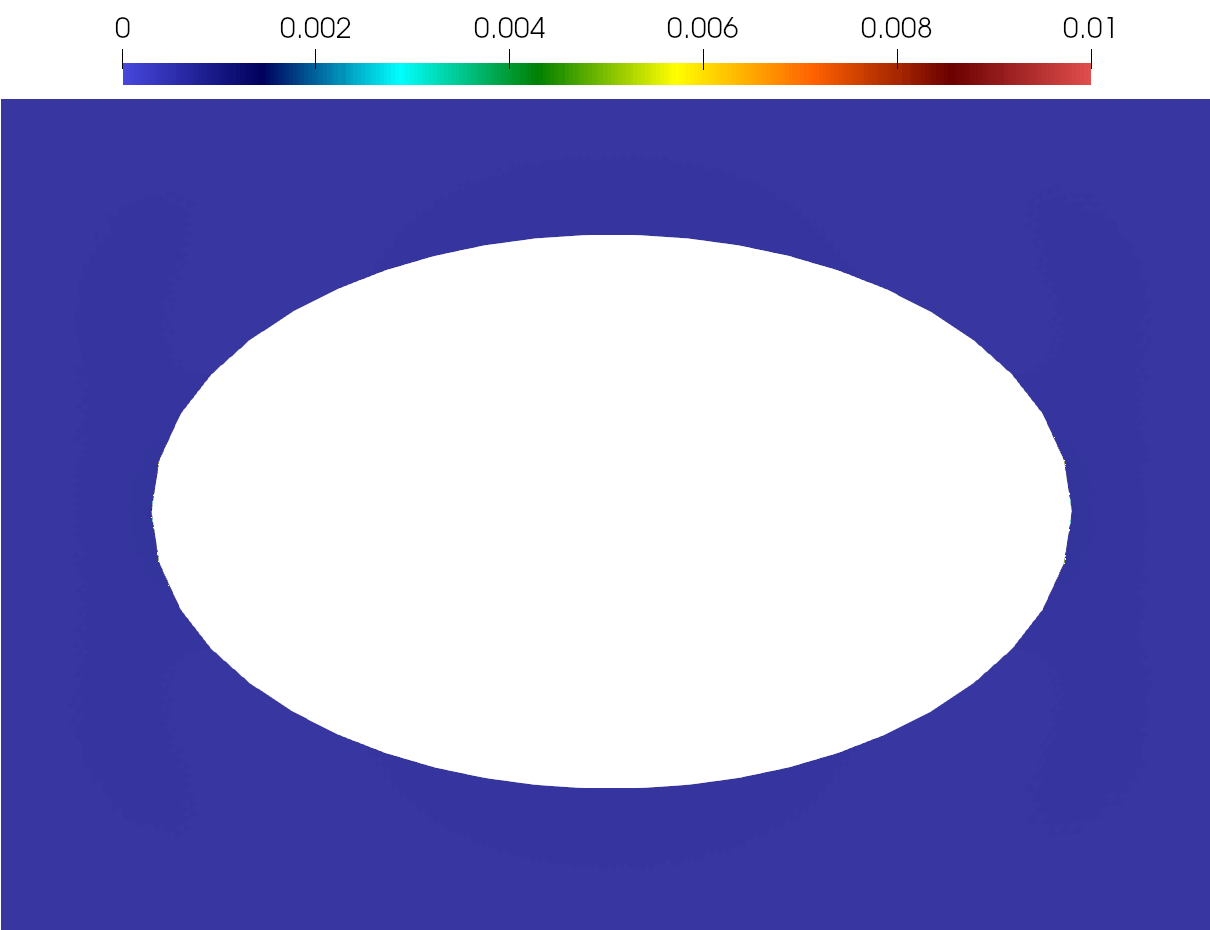}%
	}}
      \caption{Magnitude of $\linstrain$ at the final time near the interior hole.}%
      \label{fig:ellipse_residual_strain}
\end{figure}

Figure \ref{fig:ellipse_stress} shows plots of the magnitude of the deviatoric part of the stress $\linstress_\delta$, for the elastic problem and the problem with $\kappa_\star$; the domain is deformed according to the solution for the displacement. (The deformation is magnified 15 times.) We observe stresses concentrating on the sides of the elliptical hole. Note that the stress reaches much higher values in the elastic case. Figure \ref{fig:ellipse_residual_strain} shows the solution at the final time $t=1$\,s; in the elastic case the strain is practically zero, while in the plastic case there is still a residual strain on the sides of the elliptical hole.

\begin{figure}
\centering
\subfloat[C][{\centering $t=0$\,s}]{{%
	\includegraphics[width=0.42\textwidth]{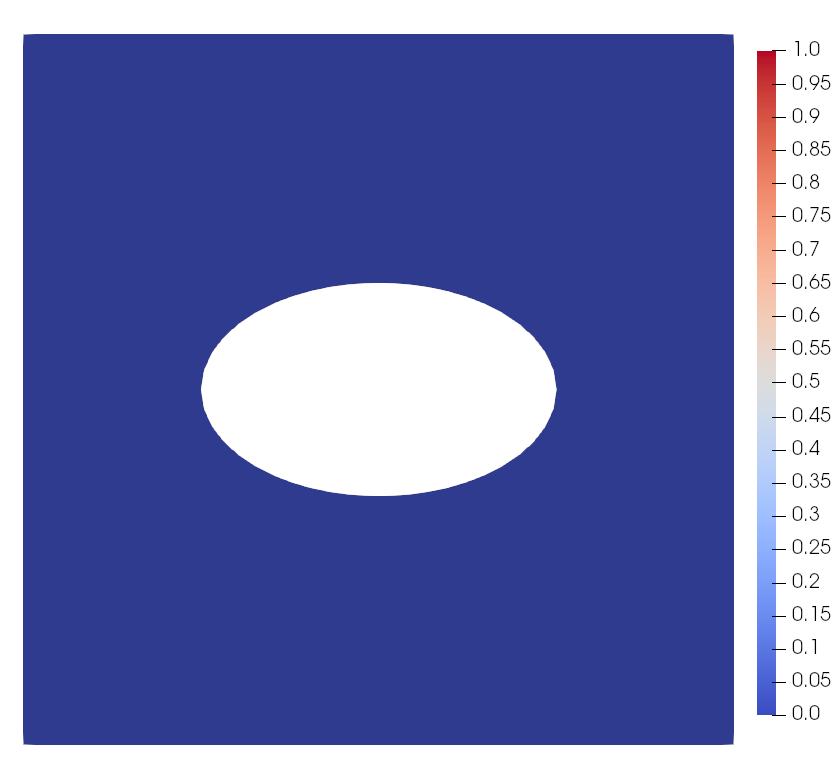}%
	}}%
	\subfloat[D][{\centering $t=0$\,s}]{{%
	\includegraphics[width=0.42\textwidth]{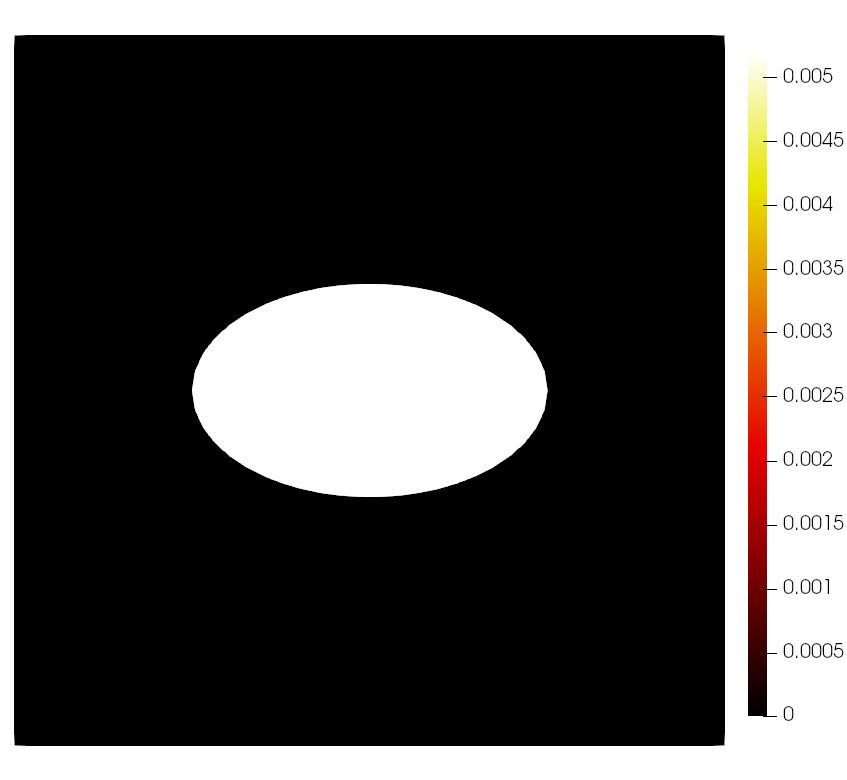}%
	}}\\
\subfloat[C][{\centering $t=0.375$\,s}]{{%
	\includegraphics[width=0.42\textwidth]{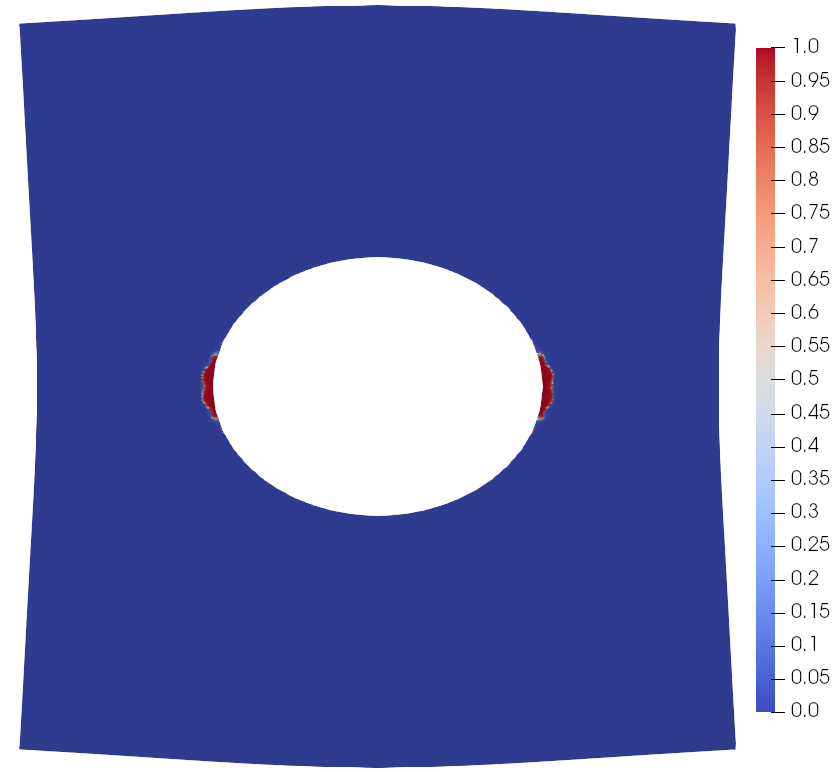}%
	}}%
	\subfloat[D][{\centering $t=0.375$\,s}]{{%
	\includegraphics[width=0.42\textwidth]{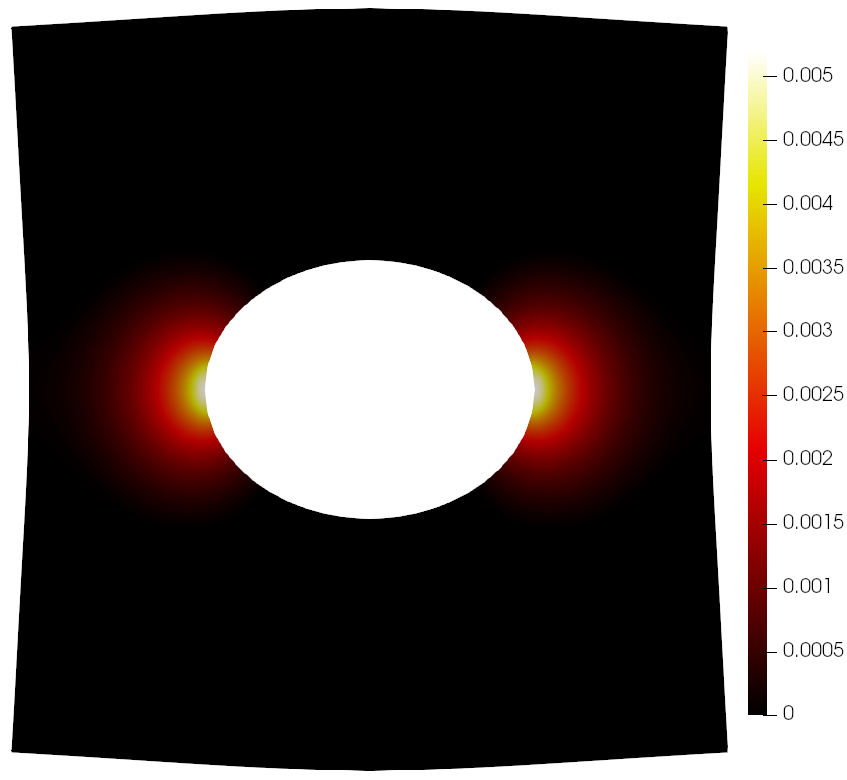}%
	}}\\
\subfloat[C][{\centering $t=0.54$\,s}]{{%
	\includegraphics[width=0.42\textwidth]{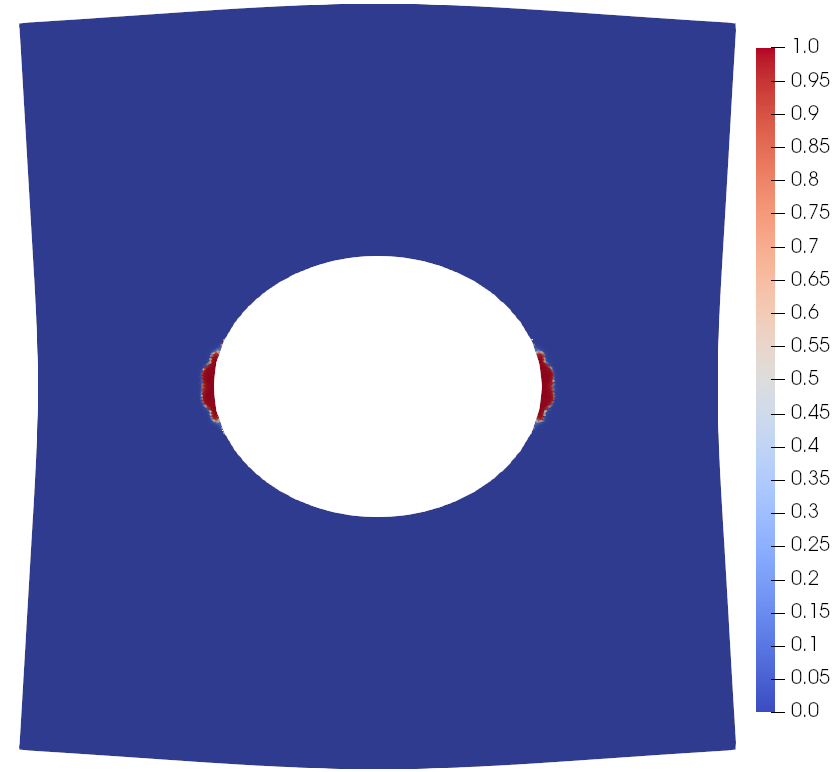}%
	}}%
	\subfloat[D][{\centering $t=0.54$\,s}]{{%
	\includegraphics[width=0.42\textwidth]{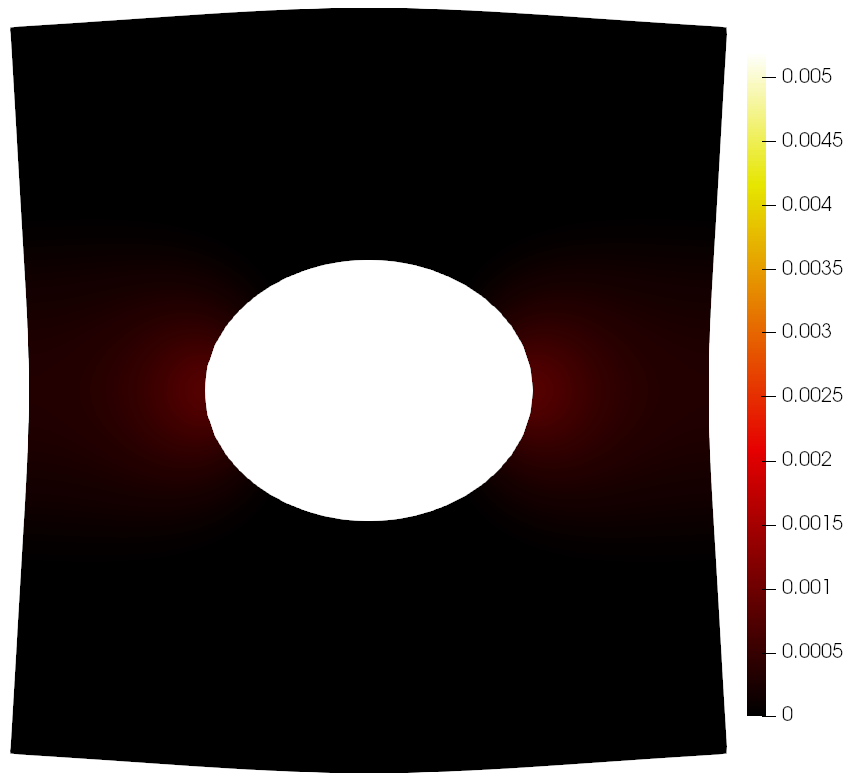}%
	}}\\
	\caption{Plot of the temperature difference $\theta$ with respect to the initial state (right) and $H_\epsilon(|\linstress_\delta|^2-\kappa_\star^2)$ (left) for the problem with $\kappa_\star=60$\,Pa.}%
	\label{fig:ellipse_temperature}
\end{figure}

The evolution of the temperature difference $\theta$ with respect to the initial temperature is shown in Figure~\ref{fig:ellipse_temperature}, along with plots of the function $H_\epsilon(|\linstress_\delta|^2 - \kappa_\star^2)$, which allows us to track whether the yield criterion is satisfied. In Figure~\ref{fig:ellipse_temperature} we observe the behaviour expected from the model; namely, the regions where the stress concentrates---and where the yield criterion is satisfied---act as a heat source for the temperature field. Note that without this heat source the temperature field would be otherwise identically zero, thanks to the boundary conditions. Moreover, we also see in Figure~\ref{fig:ellipse_temperature} (\textsc{F}) that at time $t=0.54$\,s there is no longer a heat source for the temperature field, since at this time the loading criterion $\linstress\fp \linstrain >0$ is not satisfied.

\section{Concluding remarks}
The family of models stemming from the simple rate-type equation~\eqref{eq:1} provides a novel approach to the modelling of inelastic response. Starting with the introduction of the simple rate-type equation~\eqref{eq:1} in~\cite{rajagopal.kr.srinivasa.ar:inelastic}, several variants and generalisation of~\eqref{eq:1} have been investigated, see, for example, \cite{rajagopal.kr.srinivasa.ar:implicit*1}, and a thermodynamical framework for some of these models has been successfully developed even in the finite deformations setting, see~\cite{CP.2020}. (For another treatment of rate-type models from a thermodynamic point of view see also~\cite{giorgi.c.morro.a:thermodynamic} and the discussion in~\cite{houlsby.gt.puzrin.am:principles} based on the concept of internal variables, to name a few.) Furthermore, the models in this class have also been employed in modelling the response of various materials, see~\cite{wang.z.srinivasa.ar.ea:simulation}, \cite{mozafari.f.thamburaja.p.ea:rate}, \cite{mozafari.f.thamburaja.p.ea:finite-element} \cite{saravanan.u.rajagopal.kr.ea:model} and~\cite{bustamante.r.rajagopal.kr:three-dimensional}. We have focused on a generalisation of~\eqref{eq:1} that describes the standard elastic--perfectly plastic response.

In particular, we have investigated numerical schemes for the solution of the corresponding governing equations~\eqref{eq:pde_Tv} and~\eqref{eq:pde_ut} respectively. Given the novelty of the model, the mathematical theory for the corresponding model is clearly underdeveloped compared to the mathematical theory for the classical models of elastic--perfectly plastic behaviour, see, for example, \cite{DL.1976}, \cite{SH.2006} or~\cite{kruzk.m.roubcek.t:mathematical}. From this perspective, it might seem useless to develop yet another variant of mathematical theory for the standard elastic--perfectly plastic response. However, our analysis serves a different purpose. We focus on a prototypical example of a rate-type evolution equation for a rate-independent process, and our objective is to investigate the viability of the rate-type models based approach. Naturally, the vision is to continue with numerical analysis of more involved models for inelastic responses that go beyond the standard elastic--perfectly plastic response.  

The considered model for elastic--perfectly plastic response is from the physical point of view conceptually very clean and simple, and we show in this article that this transfers to the numerical analysis as well. The model is amenable to standard discretisation techniques, and we show that the ``straightforward'' finite element discretisation of the model inherits directly the energy stability properties of the continuous model. Furthermore, since the model possesses a solid thermodynamical basis, we can also formulate the corresponding temperature evolution equation and provide numerical analysis for the full thermo-mechanical problem. (Concerning the temperature evolution, our model is \textcolor{black}{however} a simple one, the thermal response related to the \emph{plastic} deformation can be more complicated, see \cite{rosakis.p.rosakis.aj.ea:thermodynamic} and subsequent works; the same holds for the \emph{elastic} response, complexities such as Gough--Joule effect, see \cite{gough.j:description}, \cite{joule.jp:on*2} and modern treatment thereof in~\cite{anand.l:constitutive}, are in the present rigorous numerical analysis neglected. On the level of mathematical modelling, it is possible to formulate thermodynamically consistent models with more involved thermal effects such as models with temperature dependent elastic/plastic parameters, but for such models the numerical analysis would be more demanding.) Finally, we show that---up to numerical dissipation---all the energy budget of the system is accounted for. The numerical analysis is documented by an implementation of the proposed scheme.

One of the main takeaways of this work is that the rate-type model under consideration is numerically tractable via standard numerical schemes without the need for more sophisticated methods required for instance by variational inequality formulations. From a broader perspective this suggests that the modelling of inelastic rate-independent phenomena via the rate-type models might be---from the theoretical numerical analysis point of view---feasible as well. In particular numerical schemes for models describing complex inelastic phenomena such as the Mullins effect, see~\cite{diani.j.fayolle.b.ea:review} for a general discussion and~\cite{cichra.d.gazca-orozco.pa.ea:thermodynamic} for a rate-type model, might be of interest in this regard. Our work is a precursor for such studies.

\bibliographystyle{wileyNJD-MPS}
\bibliography{./vit-prusa,./bibliography}
\end{document}